\documentclass[thmsa,final,11pt]{article}

\usepackage{eurosym}
\usepackage{amsfonts}
\usepackage{amsmath}
\usepackage{amssymb}
\usepackage{mathrsfs} 
\usepackage{enumerate}
\usepackage{graphicx}
\usepackage{subfigure}
\usepackage{color,xcolor}
\usepackage[onehalfspacing,nodisplayskipstretch]{setspace}
\usepackage[colorlinks=true, linkcolor=blue, citecolor=blue, urlcolor=blue]{hyperref}
\numberwithin{equation}{section}

\newcommand{\hide}[1]{}

\newtheorem{theorem}{Theorem}

\newtheorem{lemma}[theorem]{Lemma}

\newtheorem{proposition}[theorem]{Proposition}

\normalsize
\newenvironment{prooft}[1]{\vskip 2mm\noindent {\bf Proof of #1.}}
                    {\hfill $\square$ \vskip 2mm \noindent}
\usepackage[toc,page]{appendix}

\begin{document}

\title{Central limit theorems for squared increment sums of fractional Brownian fields based on a Delaunay triangulation in $2D$}
\author{Nicolas CHENAVIER\thanks{%
Universit\'{e} du Littoral C\^{o}te d'Opale, 50 rue F. Buisson 62228 Calais.
nicolas.chenavier@univ-littoral.fr}\quad and\quad Christian Y.\ ROBERT%
\thanks{%
1. Universit\'{e} de
Lyon, Universit\'{e} Lyon 1, Institut de Science Financi\`{e}re et
d'Assurances, 50 Avenue Tony Garnier, F-69007 Lyon, France.
2. Laboratory in Finance and Insurance - LFA CREST - Center for Research in
Economics and Statistics, ENSAE, Palaiseau, France. christian.robert@univ-lyon1.fr%
}}
\maketitle

\abstract{We study quadratic variations of a two-dimensional isotropic fractional Brownian field observed on a random spatial design. The observation sites are given by a homogeneous Poisson point process with intensity \(N\) in a fixed unit square, and increments are computed along the edges of the associated Poisson--Delaunay triangulation. This random geometric framework differs from the regular-grid settings usually considered for quadratic variations of Gaussian fields.
For Hurst parameter \(H<1/2\), we establish central limit theorems, as \(N\to\infty\), for two classes of centered squared-increment statistics. The first one is based on normalized increments along Delaunay edges, whereas the second one is based on normalized pairs of increments associated with two edges of each Delaunay triangle. The asymptotic variances are finite and admit integral representations involving the geometry of the typical Poisson--Delaunay tessellation.
The proofs combine the scaling properties of fractional Brownian fields, ergodic arguments for marked Poisson point processes, and a Malliavin--Stein version of the Breuer--Major theorem. A simulation study illustrates the Gaussian approximation in the admissible range \(H<1/2\) and the emergence of non-Gaussian behavior beyond the critical threshold.

}


\strut

\textit{AMS (2020):} 60F05, 60G22, 60D05, 60G55.

\maketitle

\section{Introduction}

Fractional Brownian motion is a centered Gaussian self-similar process with
stationary increments. It provides a canonical model for phenomena exhibiting
rough sample paths and medium- or long-range dependence, and has therefore
been widely used in applications where classical Brownian, semimartingale or
Markovian models are too restrictive; examples include financial data and
network traffic, where persistence effects play an important role
(see e.g. \cite{Mikosch2002}).

There is, however, no unique extension of fractional Brownian motion to
random fields indexed by multidimensional spaces. The reason is that, for
random fields, several notions of increments may be considered
(see e.g. Section~3.3 in \cite{Cohen&Istas13}). In this paper we focus on
the two-dimensional case. A first natural notion is that of linear stationary
increments. A random field
\(W:=\left(W(x)\right)_{x\in\mathbf{R}^{2}}\) is said to have linear
stationary increments if the law of
\[
 \left(W(x+x_{0})-W(x_{0})\right)_{x\in\mathbf{R}^{2}}
\]
does not depend on \(x_{0}\in\mathbf{R}^{2}\). The main example considered
throughout this paper is the isotropic fractional Brownian field, defined as
the centered Gaussian random field satisfying \(W(0)=0\) a.s. and
\begin{equation}
\mathrm{cov}(W(x),W(y))
 =
 \frac{\sigma^{2}}{2}
 \left(
   \|x\|^{2H}+\|y\|^{2H}-\|y-x\|^{2H}
 \right),
 \label{eq:defcovariance}
\end{equation}
where \(H\in(0,1)\), \(\sigma^{2}>0\), and \(\|\cdot\|\) denotes the
Euclidean norm. The parameter \(\sigma\) is the scale parameter, whereas
\(H\) is the Hurst parameter and governs the local regularity of the field.

Another classical notion is based on rectangular increments. For
\(x=(x_{1},x_{2})\in\mathbf{R}^{2}\), \(x_{0}\in\mathbf{R}^{2}\),
\(e_{1}=(1,0)\) and \(e_{2}=(0,1)\), the rectangular increment of \(W\) is
defined by
\begin{equation}
W(x+x_{0})
 -W(x_{1}e_{1}+x_{0})
 -W(x_{2}e_{2}+x_{0})
 +W(x_{0}).
\label{eq:defrectinc}
\end{equation}
A field has rectangular stationary increments if the law of the process of
rectangular increments in \eqref{eq:defrectinc} does not depend on the
choice of \(x_{0}\). The fractional Brownian sheet is a standard example of
a self-similar Gaussian field with rectangular stationary increments. This
distinction between linear and rectangular increments is important because it
leads to different notions of quadratic variation for random fields. In the
present paper, we work with the isotropic fractional Brownian field
\eqref{eq:defcovariance} and with increments computed along line segments.

Quadratic variations, or squared increment sums, play a central role in
stochastic analysis and statistical inference. For stochastic processes, they
are fundamental both from a theoretical point of view
(see e.g. \cite{Protter2005}, pp.~66--77) and for the estimation of model
parameters. They have been extensively studied for fractional Brownian motion
(see e.g. \cite{Nourdin08,Nourdin10}) and, more generally, for Gaussian
sequences (see the survey \cite{Viitasaari19}). For random fields, the
definition of quadratic variation depends on the chosen notion of increments.
Quadratic variations based on rectangular increments over regular grids were
studied, for instance, in \cite{Deo&Wong79}, while functional limit theorems
for generalized variations of the fractional Brownian sheet were obtained in
\cite{Pakkanen&Reveillac16}. For isotropic fractional Brownian fields,
linearly filtered increments over regular grids were introduced in
\cite{Zhu&Stein02} to estimate fractal dimension.

The literature is much more limited for irregularly spaced observations.
In \cite{Loh15}, observations are taken along a smooth curve in space, and
higher-order quadratic variations are used to estimate the smoothness of
Gaussian random fields. To the best of our knowledge, no previous work has
established central limit theorems for quadratic variations of fractional
Brownian fields observed at random spatial locations and based on the
geometry induced by those locations.

The purpose of this paper is to fill this gap. We consider an isotropic
fractional Brownian field \(W\) with covariance \eqref{eq:defcovariance},
independent of a homogeneous Poisson point process \(P_{N}\) with intensity
\(N\) in \(\mathbf{R}^{2}\). The field is observed at the points of \(P_{N}\)
lying in the fixed unit square $\mathbf{C}=(-1/2,1/2]^{2}$,
and increments are computed along edges of the associated Poisson--Delaunay
triangulation, denoted by \(\operatorname{Del}(P_{N})\). Thus, unlike the
regular-grid framework, both the observation sites and the set of increments
are random. Our main contribution is to establish central limit theorems, as
\(N\to\infty\), for centered squared normalized increment sums built from
Delaunay edges and from pairs of Delaunay edges within triangles.

The use of the Delaunay triangulation is motivated by a statistical problem.
In \cite{Chenavier&Robert25c}, composite maximum likelihood estimators for
max-stable Brown--Resnick random fields are constructed from pairs and
triples of observation sites selected through the Delaunay triangulation.
This choice is natural because the Delaunay triangulation is, in a precise
geometric sense, one of the most regular triangulations: it maximizes the
minimum angle among triangulations of the same point configuration. The
asymptotic analysis of such estimators requires central limit theorems for
statistics associated with a single fractional Brownian field observed on a
random design. This is the problem addressed here. The Poisson assumption
provides a tractable infill model and allows explicit computations through
the Slivnyak--Mecke formula.

\paragraph{Squared normalized increment sums}

We now introduce the two statistics studied in the paper. Recall that
\(W\) is the isotropic fractional Brownian field with covariance
\eqref{eq:defcovariance}, and that \(P_{N}\) is a homogeneous Poisson point
process with intensity \(N\), independent of \(W\). When
\(x_{1},x_{2}\in P_{N}\) are Delaunay neighbors, we write
\(x_{1}\sim x_{2}\) in \(\operatorname{Del}(P_{N})\). Let \(E_{N}\) be the
set of ordered Delaunay edges \((x_{1},x_{2})\) such that
\[
x_{1}\sim x_{2}\text{ in }\operatorname{Del}(P_{N}),
\qquad
x_{1}\in\mathbf{C},
\qquad
x_{1}\preceq x_{2},
\]
where \(\preceq\) denotes the lexicographic order. Similarly, let \(DT_{N}\)
be the set of ordered Delaunay triangles \((x_{1},x_{2},x_{3})\) such that
\[
\Delta(x_{1},x_{2},x_{3})\in\operatorname{Del}(P_{N}),
\qquad
x_{1}\in\mathbf{C},
\qquad
x_{1}\preceq x_{2}\preceq x_{3},
\]
where \(\Delta(x_{1},x_{2},x_{3})\) denotes the convex hull of
\(\{x_{1},x_{2},x_{3}\}\).

For two distinct points \(x_{1},x_{2}\in\mathbf{R}^{2}\), write
\(d_{1,2}=\|x_{2}-x_{1}\|\), and define the normalized increment
\[
U_{x_{1},x_{2}}^{(W)}
 =
 \sigma^{-1}d_{1,2}^{-H}
 \left(W(x_{2})-W(x_{1})\right).
\]
This normalization ensures that \(U_{x_{1},x_{2}}^{(W)}\) has a standard
Gaussian distribution.

The first statistic is the centered squared increment sum over Delaunay
edges:
\begin{equation*}
V_{2,N}^{(W)}
 =
 \frac{1}{\sqrt{|E_{N}|}}
 \sum_{(x_{1},x_{2})\in E_{N}}
 \left\{
   \left(U_{x_{1},x_{2}}^{(W)}\right)^{2}-1
 \right\}.
\end{equation*}
Here \(|E_{N}|\) denotes the cardinality of \(E_{N}\).

The second statistic is based on pairs of edges in Delaunay triangles. For
\((x_{1},x_{2},x_{3})\in DT_{N}\), set
\begin{equation}
R_{x_{1},x_{2},x_{3}}
 =
 \mathrm{corr}\left(
   U_{x_{1},x_{2}}^{(W)},U_{x_{1},x_{3}}^{(W)}
 \right)
 =
 \frac{
   d_{1,2}^{2H}+d_{1,3}^{2H}-d_{2,3}^{2H}
 }{
   2(d_{1,2}d_{1,3})^{H}
 },
\label{eq:coorU}
\end{equation}
where \(d_{1,3}=\|x_{3}-x_{1}\|\) and \(d_{2,3}=\|x_{3}-x_{2}\|\). We define
\begin{multline*}
V_{3,N}^{(W)}
 =
 \frac{1}{\sqrt{|DT_{N}|}}
 \sum_{(x_{1},x_{2},x_{3})\in DT_{N}}
 \Bigg[
 \left(
 \begin{array}{cc}
 U_{x_{1},x_{2}}^{(W)}
 &
 U_{x_{1},x_{3}}^{(W)}
 \end{array}
 \right)
 \left(
 \begin{array}{cc}
 1 & R_{x_{1},x_{2},x_{3}} \\
 R_{x_{1},x_{2},x_{3}} & 1
 \end{array}
 \right)^{-1}
 \left(
 \begin{array}{c}
 U_{x_{1},x_{2}}^{(W)} \\
 U_{x_{1},x_{3}}^{(W)}
 \end{array}
 \right)
 -2
 \Bigg].
\end{multline*}
Only two edges of each triangle are used, namely \([x_{1},x_{2}]\) and
\([x_{1},x_{3}]\), since the third increment along \([x_{2},x_{3}]\) is
determined by the first two.

Equivalently, \(V_{3,N}^{(W)}\) can be written as a sum of centered squared
standard Gaussian variables. Indeed, set
\[
\widetilde{U}_{x_{1},x_{2},x_{3}}^{(W)}
 =
 (1-R_{x_{1},x_{2},x_{3}}^{2})^{-1/2}
 \left(
   U_{x_{1},x_{2}}^{(W)}
   -
   R_{x_{1},x_{2},x_{3}}U_{x_{1},x_{3}}^{(W)}
 \right),
\qquad
\widetilde{U}_{x_{1},x_{3}}^{(W)}
 =
 U_{x_{1},x_{3}}^{(W)}.
\]
Then
\[
V_{3,N}^{(W)}
 =
 \frac{1}{\sqrt{|DT_{N}|}}
 \sum_{(x_{1},x_{2},x_{3})\in DT_{N}}
 \left[
   \left(\widetilde{U}_{x_{1},x_{2},x_{3}}^{(W)}\right)^{2}-1
   +
   \left(\widetilde{U}_{x_{1},x_{3}}^{(W)}\right)^{2}-1
 \right].
\]
The construction follows the orthogonalization used in \cite{Chan&Wood00}:
for each triangle,
\[
\mathrm{corr}\left(
 \widetilde{U}_{x_{1},x_{2},x_{3}}^{(W)},
 \widetilde{U}_{x_{1},x_{3}}^{(W)}
\right)=0.
\]

Strictly speaking, the statistics above are not defined on the events
\(\{|E_N|=0\}\) and \(\{|DT_N|=0\}\). We shall use the convention that
\(V_{2,N}^{(W)}=0\) on \(\{|E_N|=0\}\) and \(V_{3,N}^{(W)}=0\) on
\(\{|DT_N|=0\}\). This convention is asymptotically immaterial, since these
exceptional events have exponentially small probability as \(N\to\infty\).

\paragraph{Main result}

The main result of the paper is the following central limit theorem.

\begin{theorem}
\label{th:CLTgaussian}
Let \(W\) be an isotropic fractional Brownian field with covariance
\eqref{eq:defcovariance}, where \(H\in(0,1/2)\) and \(\sigma^{2}>0\). There
exist finite constants \(\sigma_{V_{2}}^{2}>0\) and
\(\sigma_{V_{3}}^{2}>0\) such that, as \(N\to\infty\),
\[
V_{2,N}^{(W)}
 \overset{\mathcal{D}}{\longrightarrow}
 \mathcal{N}(0,\sigma_{V_{2}}^{2}),
\qquad
V_{3,N}^{(W)}
 \overset{\mathcal{D}}{\longrightarrow}
 \mathcal{N}(0,\sigma_{V_{3}}^{2}).
\]
\end{theorem}

The normalization by the square root of the number of Delaunay edges or
triangles is the same as in the regular-grid central limit theorems of
\cite{Chan&Wood00} and \cite{Zhu&Stein02}. The asymptotic variances in
Theorem~\ref{th:CLTgaussian} are finite and admit explicit integral
representations involving the geometry of the Poisson--Delaunay
triangulation, although these expressions are rather involved.

The theorem is stated in a fixed-domain, or infill, asymptotic framework:
the observation window is fixed, while the intensity \(N\) of the Poisson
point process tends to infinity. By scaling, this is equivalent to considering
a Poisson point process with fixed intensity on an expanding window. We use
the fixed-domain formulation because it is the relevant one for the
companion papers \cite{Chenavier&Robert25b,Chenavier&Robert25c}.

The proof combines three main ingredients. First, the scaling property of
the fractional Brownian field reduces the problem to an expanding-window
setting with a unit-intensity Poisson point process. Second, ergodic
arguments for suitable marked Poisson point processes provide the
limits of the conditional variances. Third, a Malliavin--Stein version of the
Breuer--Major theorem, due to Nourdin and Peccati
(Theorem~7.2.4 in \cite{Nourdin&Peccati12}), yields the Gaussian limit once
the relevant chaos contraction conditions have been checked. The main
technical difficulty is to control correlations between normalized increments
indexed by a random geometric graph, including boundary effects. The theorem
is stated for edges and triangles whose lexicographically first vertex lies
in \(\mathbf{C}\); analogous results can be obtained for configurations
fully contained in \(\mathbf{C}\) by treating the corresponding edge effects.

The restriction \(H<1/2\) is intrinsic to the order of the increments used in
this paper. For increments of order \(0\), this threshold is consistent
with the known regime change for quadratic variations of fractional fields
on regular grids; see, for example, Remark~3.2 in \cite{Chan&Wood00}. In the
univariate case, non-central limits of Rosenblatt type arise under suitable
renormalization in the long-memory regime \cite{Dobrushin&Major79}. It is
therefore natural to expect a non-Gaussian limit for the Delaunay-based
statistics considered here when \(H>1/2\). By contrast, using increments of
order \(1\), which would require configurations involving four vertices, is
expected to restore Gaussian limits for all \(H\in(0,1)\), as happens in the
regular-grid setting \cite{Bierme11,Chan&Wood04}.
The boundary case \(H=1/2\) is not covered by Theorem~1. It corresponds to a
critical regime for the squared-increment statistics considered here. Indeed,
the correlations between distant normalized linear increments decay as
\(\|x-y\|^{2H-2}\), so that the covariances of squared increments decay as
\(\|x-y\|^{4H-4}\). In dimension two, this decay is integrable at infinity
only when \(H<1/2\). At \(H=1/2\), the corresponding integral diverges
logarithmically. One may therefore expect a Gaussian limit to persist at the
critical value, but with an additional logarithmic normalization, namely with
\(\sqrt{|E_N|\log N}\) and \(\sqrt{|DT_N|\log N}\) instead of
\(\sqrt{|E_N|}\) and \(\sqrt{|DT_N|}\). Establishing such a critical central
limit theorem would require a separate analysis of the logarithmically
divergent variance and is left for future work.

The paper is organized as follows. Section~\ref{sec:preliminaries} recalls
the basic facts on Poisson--Delaunay triangulations used throughout the
paper. Section~\ref{sec:proofedges} gives the detailed proof of
Theorem~\ref{th:CLTgaussian} for the edge-based statistic
\(V_{2,N}^{(W)}\). Section~\ref{sec:simulations} presents a Monte Carlo
study illustrating the finite-sample behavior of the statistic and the
sharpness of the condition \(H<1/2\). Appendix~\ref{sec:proofpairsedges}
sketches the proof for the triangle-based statistic \(V_{3,N}^{(W)}\), whose
arguments are similar. Appendix~\ref{sec:intermediary_results} collects
technical estimates used in the proof.

\section{Preliminaries}
\label{sec:preliminaries}

This section recalls the basic facts on Poisson--Delaunay triangulations
that will be used throughout the paper. We follow the terminology and
normalizations of \cite{Schneider&Weil08}.

Let \(P_N\) be a homogeneous Poisson point process with intensity \(N\) in
\(\mathbf{R}^{2}\). The Delaunay triangulation associated with \(P_N\),
denoted by \(\operatorname{Del}(P_N)\), is the triangulation with vertex set
\(P_N\) such that the circumdisk of each triangle contains no point of
\(P_N\) in its interior. For a homogeneous Poisson point process, this
triangulation is well defined and unique almost surely; see, for instance,
p.~478 in \cite{Schneider&Weil08}. The Delaunay triangulation is a canonical
choice in computational geometry because it avoids thin triangles as much as
possible: among triangulations of a fixed point configuration, it maximizes
the minimum angle.

We first recall the notion of typical cell for the Poisson--Delaunay
triangulation generated by a unit-intensity Poisson point process \(P_1\).
For each triangular cell \(C\in \operatorname{Del}(P_1)\), let \(z(C)\)
denote its circumcenter. If \(\mathbf{B}\subset \mathbf{R}^{2}\) is a Borel
set with area \(a(\mathbf{B})\in(0,\infty)\), the cell intensity
\(\beta_2\) is defined by
\[
\beta_2
 =
 \frac{1}{a(\mathbf{B})}
 \mathbb{E}
 \left[
   \left|
     \{C\in \operatorname{Del}(P_1): z(C)\in \mathbf{B}\}
   \right|
 \right].
\]
It is well known that \(\beta_2=2\); see Theorem~10.2.9 in
\cite{Schneider&Weil08}.

The typical cell is a random triangle \(\mathcal{C}\) whose distribution is
defined by the Palm-type identity
\[
\mathbb{E}\left[g(\mathcal{C})\right]
 =
 \frac{1}{\beta_2 a(\mathbf{B})}
 \mathbb{E}
 \left[
   \sum_{C\in \operatorname{Del}(P_1): z(C)\in \mathbf{B}}
   g(C)
 \right],
\]
for every positive measurable translation-invariant function
\(g:\mathcal{K}_2\to\mathbf{R}\), where \(\mathcal{K}_2\) denotes the space
of compact convex subsets of \(\mathbf{R}^{2}\), endowed with the Fell
topology. This definition does not depend on the particular choice of
\(\mathbf{B}\). The distribution of \(\mathcal{C}\) admits the following
integral representation; see Theorem~10.4.4 in \cite{Schneider&Weil08}:
\begin{equation}
\mathbb{E}\left[g(\mathcal{C})\right]
 =
 \frac{1}{6}
 \int_{0}^{\infty}
 \int_{(\mathbf{S}^{1})^{3}}
 r^{3}e^{-\pi r^{2}}
 a\!\left(\Delta(u_{1},u_{2},u_{3})\right)
 g\!\left(\Delta(ru_{1},ru_{2},ru_{3})\right)
 \sigma(\mathrm{d}u_{1})
 \sigma(\mathrm{d}u_{2})
 \sigma(\mathrm{d}u_{3})
 \mathrm{d}r .
\label{eq:typicalcell}
\end{equation}
Here \(\mathbf{S}^{1}\) is the unit circle of \(\mathbf{R}^{2}\), and
\(\sigma\) is the spherical Lebesgue measure normalized by
\(\sigma(\mathbf{S}^{1})=2\pi\). Equivalently,
\(\mathcal{C}\) has the same distribution as
\(R\Delta(U_{1},U_{2},U_{3})\), where \(R\) is independent of
\((U_{1},U_{2},U_{3})\), with respective densities
\[
2\pi^{2}r^{3}e^{-\pi r^{2}},
\qquad r>0,
\]
and
\[
\frac{a(\Delta(u_{1},u_{2},u_{3}))}{12\pi^{2}},
\qquad
(u_{1},u_{2},u_{3})\in(\mathbf{S}^{1})^{3}.
\]

We shall also use the corresponding notion of typical edge. The edge
intensity \(\beta_1\) of \(\operatorname{Del}(P_1)\) is the mean number of
edges per unit area. It is equal to \(\beta_1=3\); see again
Theorem~10.2.9 in \cite{Schneider&Weil08}. If \(D\) denotes the length of
the typical edge, then \(D\) has the same distribution as
\(R\|U_{1}-U_{2}\|\). Its distribution function satisfies, for every
\(\ell>0\),
\begin{align}
\mathbb{P}(D\leq \ell)
&=
\int_{0}^{\ell} f_D(d)\,\mathrm{d}d
\notag\\
&=
\frac{\pi}{3}
\int_{0}^{\infty}
\int_{(\mathbf{S}^{1})^{2}}
r^{3}e^{-\pi r^{2}}
a\!\left(\Delta(u_{1},u_{2},e_{1})\right)
\mathbb{I}
\left[
  r\|u_{1}-u_{2}\|\leq \ell
\right]
\sigma(\mathrm{d}u_{1})
\sigma(\mathrm{d}u_{2})
\mathrm{d}r ,
\label{eq:typicallength}
\end{align}
where \(e_{1}=(1,0)\).

In the analysis of the statistic based on pairs of edges, we also need a
typical pair of distinct Delaunay edges sharing one common vertex. Following
\eqref{eq:typicalcell}, this object may be represented by a random vector
\((D_{1},D_{2},\Theta)\), where \(D_{1},D_{2}\geq0\) are the lengths of the
two edges and \(\Theta\in[-\pi/2,\pi/2)\) is the angle between them. Its
distribution is given by
\begin{multline*}
\mathbb{P}\left[(D_{1},D_{2},\Theta)\in B\right]
 =
 \frac{1}{6}
 \int_{0}^{\infty}
 \int_{(\mathbf{S}^{1})^{3}}
 r^{3}e^{-\pi r^{2}}
 a\!\left(\Delta(u_{1},u_{2},u_{3})\right)
\\
\times
\mathbb{I}
\left[
 \left(
   r\|u_{3}-u_{2}\|,
   r\|u_{2}-u_{1}\|,
   \arcsin\left(\cos(\theta_{u_{1},u_{2}}/2)\right)
 \right)
 \in B
\right]
\sigma(\mathrm{d}u_{1})
\sigma(\mathrm{d}u_{2})
\sigma(\mathrm{d}u_{3})
\mathrm{d}r ,
\end{multline*}
for every Borel set
\(B\subset \mathbf{R}_{+}^{2}\times[-\pi/2,\pi/2)\), where
\(\theta_{u_{1},u_{2}}\) denotes the angle between \(u_{1}\) and \(u_{2}\).
Thus \(D_1\) and \(D_2\) represent the lengths of the two typical adjacent
edges, whereas \(\Theta\) represents their angle.

Throughout the paper, we identify \(\operatorname{Del}(P_N)\) with its
skeleton. Thus, when \(x_{1},x_{2}\in P_N\) are Delaunay neighbors, we write
\(x_{1}\sim x_{2}\) in \(\operatorname{Del}(P_N)\). For a Borel set
\(\mathbf{B}\subset\mathbf{R}^{2}\), let \(E_{N,\mathbf{B}}\) be the set of
ordered Delaunay edges \((x_{1},x_{2})\) satisfying
\[
x_{1}\sim x_{2}\text{ in }\operatorname{Del}(P_N),
\qquad
x_{1}\in\mathbf{B},
\qquad
x_{1}\preceq x_{2},
\]
where \(\preceq\) denotes the lexicographic order. In particular, for
\(\mathbf{C}=(-1/2,1/2]^2\), we have \(E_{N,\mathbf{C}}=E_N\).

Similarly, for a Borel set \(\mathbf{B}\subset\mathbf{R}^{2}\), let
\(DT_{N,\mathbf{B}}\) be the set of ordered Delaunay triangles
\((x_{1},x_{2},x_{3})\) satisfying
\[
\Delta(x_{1},x_{2},x_{3})\in \operatorname{Del}(P_N),
\qquad
x_{1}\in\mathbf{B},
\qquad
x_{1}\preceq x_{2}\preceq x_{3},
\]
where \(\Delta(x_{1},x_{2},x_{3})\) denotes the convex hull of
\(\{x_{1},x_{2},x_{3}\}\). When \(\mathbf{B}=\mathbf{C}\), we simply write
\(DT_{N,\mathbf{C}}=DT_N\).

\section{Proof of Theorem \protect\ref{th:CLTgaussian} for \(V_{2,N}^{(W)}\)}
\label{sec:proofedges}

We prove in this section the central limit theorem for the edge-based
statistic \(V_{2,N}^{(W)}\). 

Recall that
\[
V_{2,N}^{(W)}
=
\frac{1}{\sqrt{|E_N|}}
\sum_{(x_1,x_2)\in E_N}
\left[
  \left(U_{x_1,x_2}^{(W)}\right)^2-1
\right].
\]
The first step is to use the scaling properties of the two random objects
involved in the definition of \(V_{2,N}^{(W)}\). Since \(P_N\) is a
homogeneous Poisson point process with intensity \(N\), the rescaled point
process \(N^{1/2}P_N\) has the same distribution as \(P_1\), a homogeneous
Poisson point process with intensity \(1\). Moreover, by the self-similarity
of the fractional Brownian field,
\[
\left(W(x)\right)_{x\in\mathbf{R}^{2}}
\overset{\mathcal{D}}{=}
\left(
  N^{H/2}W(N^{-1/2}x)
\right)_{x\in\mathbf{R}^{2}}.
\]
Because the increments are normalized by the factor
\(\|x_2-x_1\|^{-H}\), this rescaling leaves the distribution of the
normalized increments unchanged. Hence \(V_{2,N}^{(W)}\) has the same
distribution as
\[
V_{2,N}^{(W)\prime}
=
\frac{1}{\sqrt{|E_N'|}}
\sum_{(x_1,x_2)\in E_N'}
\left[
  \left(U_{x_1,x_2}^{(W)}\right)^2-1
\right],
\]
where
\[
E_N'=E_{1,\mathbf{C}_N},
\qquad
\mathbf{C}_N=(-N^{1/2}/2,N^{1/2}/2]^2.
\]
Thus \(E_N'\) is the set of ordered Delaunay edges of
\(\operatorname{Del}(P_1)\) whose lexicographically first endpoint belongs
to \(\mathbf{C}_N\).

We shall prove the central limit theorem conditionally on the Poisson point
process. More precisely, we prove that, for almost every realization of
\(P_1\),
\[
V_{2,N}^{(W)\prime}
\ \Big|\ P_1
\overset{\mathcal{D}}{\longrightarrow}
\mathcal{N}(0,\sigma_{V_2}^{2}),
\]
where the limiting variance \(\sigma_{V_2}^{2}\) is deterministic. Equivalently,
for every \(u\in\mathbf{R}\),
\[
\mathbb{E}
\left[
  \exp\left(iuV_{2,N}^{(W)\prime}\right)
  \,\middle|\, P_1
\right]
\longrightarrow
\exp\left(-\frac{1}{2}\sigma_{V_2}^{2}u^2\right)
\]
for almost every realization of \(P_1\). The unconditional convergence then
follows by dominated convergence.

We now fix a realization of \(P_1\) belonging to a full-probability event on
which the Poisson--Delaunay graph is locally finite and
\[
\frac{|E_N'|}{N}\longrightarrow 3 .
\]
The latter convergence follows from the ergodic theorem for stationary
tessellations and from the fact that the edge intensity of the
Poisson--Delaunay triangulation is equal to \(3\).

Let \(E_{1,\mathbf{R}^{2}}\) denote the set of ordered Delaunay edges
\((x_1,x_2)\) of \(\operatorname{Del}(P_1)\), with \(x_1\preceq x_2\).
This set is countable. We choose a bijection
\[
\varphi:\mathbf{Z}\longrightarrow E_{1,\mathbf{R}^{2}},
\]
and write
\[
U^{(k)}=U_{\varphi(k)}^{(W)},\qquad k\in\mathbf{Z}.
\]
We also define
\[
e_N'=\{k\in\mathbf{Z}:\varphi(k)\in E_N'\}.
\]
Then \(|e_N'|=|E_N'|\), and
\[
V_{2,N}^{(W)\prime}
=
\frac{1}{\sqrt{|e_N'|}}
\sum_{k\in e_N'}
\left[
  \left(U^{(k)}\right)^2-1
\right].
\]

Conditionally on \(P_1\), the family \((U^{(k)})_{k\in\mathbf{Z}}\) is a
centered Gaussian family with unit variances. Therefore, there exist a real
separable Hilbert space \(\mathfrak{H}\), an isonormal Gaussian process
\(\{X(h):h\in\mathfrak{H}\}\), and vectors
\((\varepsilon_k)_{k\in\mathbf{Z}}\subset\mathfrak{H}\) such that
\[
U^{(k)}=X(\varepsilon_k),
\qquad
\|\varepsilon_k\|_{\mathfrak{H}}=1,
\]
and
\[
\langle\varepsilon_k,\varepsilon_l\rangle_{\mathfrak{H}}
=
\mathrm{corr}(U^{(k)},U^{(l)}),
\qquad k,l\in\mathbf{Z}.
\]
Define
\[
f_{N,2}
=
\frac{1}{\sqrt{|e_N'|}}
\sum_{k\in e_N'}
\varepsilon_k^{\otimes 2}
\in\mathfrak{H}^{\odot 2}.
\]
Since \(X(\varepsilon_k)^2-1=I_2(\varepsilon_k^{\otimes 2})\), where
\(I_2\) denotes the double Wiener--Itô integral, we have
\[
V_{2,N}^{(W)\prime}=I_2(f_{N,2}).
\]

We shall use the standard normal approximation criterion on the second
Wiener chaos. If \(f_N\in\mathfrak{H}^{\odot 2}\) satisfies
\[
2\|f_N\|_{\mathfrak{H}^{\otimes 2}}^2
\longrightarrow
\sigma^2
\]
and
\[
\|f_N\otimes_1 f_N\|_{\mathfrak{H}^{\otimes 2}}
\longrightarrow 0,
\]
then
\[
I_2(f_N)
\overset{\mathcal{D}}{\longrightarrow}
\mathcal{N}(0,\sigma^2).
\]
This is a direct consequence of the Fourth Moment Theorem, or equivalently
of the Malliavin--Stein criterion of Nourdin and Peccati (see Theorem 6.3.1 in \cite{Nourdin&Peccati12}).

It remains to verify the two conditions above for \(f_{N,2}\).
First,
\[
2\|f_{N,2}\|_{\mathfrak{H}^{\otimes 2}}^2
=
\mathbb{E}
\left[
  \left(V_{2,N}^{(W)\prime}\right)^2
  \,\middle|\, P_1
\right].
\]
Proposition~\ref{prop:variance} shows that, for almost every realization of
\(P_1\),
\[
\mathbb{E}
\left[
  \left(V_{2,N}^{(W)\prime}\right)^2
  \,\middle|\, P_1
\right]
\longrightarrow
\sigma_{V_2}^{2},
\]
where \(\sigma_{V_2}^{2}\in(0,\infty)\). The finiteness of this variance is
proved in Lemma~\ref{Le:finitevariances}.
Second, Section~\ref{Sect_Cond_c} proves that
\[
\|f_{N,2}\otimes_1 f_{N,2}\|_{\mathfrak{H}^{\otimes 2}}
\longrightarrow 0
\]
for almost every realization of \(P_1\).


\subsection{Conditional variance and asymptotic variance}
\label{sec:conditional_variance}

We now identify the limit of the conditional variance of
\(V_{2,N}^{(W)\prime}\) given the Poisson--Delaunay triangulation. This is
the variance condition required in the second-chaos normal approximation
criterion used in the proof of Theorem~\ref{th:CLTgaussian}.

Conditionally on \(P_1\), the family
$(
U_{x_1,x_2}^{(W)}
)_{(x_1,x_2)\in E_N'}$
is a centered Gaussian family with unit variances. Hence, for any two
ordered edges \((x_1,x_2)\) and \((x_3,x_4)\),
\[
\mathrm{cov}
\left(
\left(U_{x_1,x_2}^{(W)}\right)^2-1,
\left(U_{x_3,x_4}^{(W)}\right)^2-1
\right)
=
2\,
\mathrm{corr}
\left(
U_{x_1,x_2}^{(W)},U_{x_3,x_4}^{(W)}
\right)^2.
\]
It follows that
\begin{equation}
\mathbb{E}
\left[
\left(V_{2,N}^{(W)\prime}\right)^2
\,\middle|\, P_1
\right]
=
\frac{2}{|E_N'|}
\sum_{(x_1,x_2)\in E_N'}
\sum_{(x_3,x_4)\in E_N'}
\rho(x_1,x_2;x_3,x_4)^2,
\label{eq:conditional-variance-V2}
\end{equation}
where, for shortness,
\[
\rho(x_1,x_2;x_3,x_4)
=
\mathrm{corr}
\left(
U_{x_1,x_2}^{(W)},U_{x_3,x_4}^{(W)}
\right).
\]

The diagonal part of the double sum in
\eqref{eq:conditional-variance-V2} contributes exactly \(2\), since
\(\rho(x_1,x_2;x_1,x_2)=1\). The remaining terms are divided into two
classes: pairs of distinct edges sharing one endpoint and pairs of disjoint
edges. Since \(E_N'\) contains only edges whose lexicographically first
endpoint belongs to \(\mathbf{C}_N\), it is useful to first write
\begin{multline}
\mathbb{E}
\left[
\left(V_{2,N}^{(W)\prime}\right)^2
\,\middle|\, P_1
\right]
=
2
+
\frac{2}{|E_N'|}
\sum_{\substack{x_1\in \mathbf{C}_N\\
x_1\sim x_2,\ x_1\preceq x_2}}
\sum_{\substack{x_3\sim x_4\\
x_3\preceq x_4}}
\rho(x_1,x_2;x_3,x_4)^2
\mathbb{I}\left[(x_1,x_2)\neq (x_3,x_4)\right]
\\
-
\frac{2}{|E_N'|}
\sum_{\substack{x_1\in \mathbf{C}_N\\
x_1\sim x_2,\ x_1\preceq x_2}}
\sum_{\substack{x_3\notin \mathbf{C}_N\\
x_3\sim x_4,\ x_3\preceq x_4}}
\rho(x_1,x_2;x_3,x_4)^2 .
\label{eq:variance-boundary-decomposition}
\end{multline}
The last term is a boundary correction. It will be shown below that it is
negligible after normalization.

The non-diagonal term in the first line of
\eqref{eq:variance-boundary-decomposition} can be decomposed according to
the number of common endpoints of the two edges. The contribution of pairs
of distinct edges sharing one endpoint is the sum of the following four
terms:
\begin{align*}
&\frac{2}{|E_N'|}
\sum_{\substack{x_1\in\mathbf{C}_N,\ x_1\sim x_2,\ x_1\sim x_4\\
x_1\preceq x_2,\ x_1\preceq x_4}}
\rho(x_1,x_2;x_1,x_4)^2
\mathbb{I}[x_2\neq x_4]
\\
&+
\frac{2}{|E_N'|}
\sum_{\substack{x_1\in\mathbf{C}_N,\ x_1\sim x_2,\ x_2\sim x_4\\
x_1\preceq x_2,\ x_2\preceq x_4}}
\rho(x_1,x_2;x_2,x_4)^2
\mathbb{I}[x_1\neq x_4]
\\
&+
\frac{2}{|E_N'|}
\sum_{\substack{x_1\in\mathbf{C}_N,\ x_1\sim x_2,\ x_3\sim x_1\\
x_1\preceq x_2,\ x_3\preceq x_1}}
\rho(x_1,x_2;x_3,x_1)^2
\mathbb{I}[x_2\neq x_3]
\\
&+
\frac{2}{|E_N'|}
\sum_{\substack{x_1\in\mathbf{C}_N,\ x_1\sim x_2,\ x_3\sim x_2\\
x_1\preceq x_2,\ x_3\preceq x_2}}
\rho(x_1,x_2;x_3,x_2)^2
\mathbb{I}[x_1\neq x_3].
\end{align*}
The contribution of pairs of disjoint edges is
\[
\frac{2}{|E_N'|}
\sum_{\substack{x_1\in\mathbf{C}_N\\
x_1\sim x_2,\ x_1\preceq x_2}}
\sum_{\substack{x_3\sim x_4\\
x_3\preceq x_4}}
\rho(x_1,x_2;x_3,x_4)^2
\mathbb{I}
\left[
\{x_1,x_2\}\cap\{x_3,x_4\}=\emptyset
\right].
\]


We now introduce the deterministic quantities that will appear in the
limiting variance. For distinct points
\(x_1,x_2,x_3,x_4\in\mathbf{R}^2\), define
\begin{equation}
\label{eq:defp2N}
p_{2,N}(x_1,x_2,x_3,x_4)
=
\mathbb{P}
\left[
\begin{array}{l}
x_1\sim x_2,\ x_3\sim x_4
\text{ in }
\operatorname{Del}(P_N\cup\{x_1,x_2,x_3,x_4\}),\\
x_1\preceq x_2,\ x_3\preceq x_4
\end{array}
\right].
\end{equation}
This quantity is used for the case where the two edges have no common
endpoint.

To describe the cases where the two edges share one endpoint, set
\[
\mathcal{P}_2=\{(3,1),(3,2),(4,1),(4,2)\}.
\]
For \((j,i)\in\mathcal{P}_2\), the notation \(j\leftrightarrow i\) means
that the point \(x_j\) in the second edge is identified with the point
\(x_i\) in the first edge. More precisely,
\(q_{2,N}^{(j\leftrightarrow i)}\) is defined by substituting
\(x_j=x_i\) in the event of \eqref{eq:defp2N}. For example,
\[
q_{2,N}^{(3\leftrightarrow 1)}(x_1,x_2,x_4)
=
\mathbb{P}
\left[
x_1\sim x_2,\ x_1\sim x_4
\text{ in }
\operatorname{Del}(P_N\cup\{x_1,x_2,x_4\}),
\quad
x_1\preceq x_2,\ x_1\preceq x_4
\right].
\]
The other cases are defined analogously.

With this notation, the contribution of disjoint edges is described by
\[
\sigma_{0,V_2}^{2}
=
\int_{(\mathbf{R}^2)^3}
\rho(0,x_2;x_3,x_4)^2
p_{2,1}(0,x_2,x_3,x_4)
\,\mathrm{d}x_2\,\mathrm{d}x_3\,\mathrm{d}x_4.
\]
The four contributions corresponding to pairs of edges with one common
endpoint are
\begin{align*}
\sigma_{1,(3\leftrightarrow 1),V_2}^{2}
&=
\int_{(\mathbf{R}^2)^2}
\rho(0,x_2;0,x_4)^2
q_{2,1}^{(3\leftrightarrow 1)}(0,x_2,x_4)
\,\mathrm{d}x_2\,\mathrm{d}x_4,
\\
\sigma_{1,(3\leftrightarrow 2),V_2}^{2}
&=
\int_{(\mathbf{R}^2)^2}
\rho(0,x_2;x_2,x_4)^2
q_{2,1}^{(3\leftrightarrow 2)}(0,x_2,x_4)
\,\mathrm{d}x_2\,\mathrm{d}x_4,
\\
\sigma_{1,(4\leftrightarrow 1),V_2}^{2}
&=
\int_{(\mathbf{R}^2)^2}
\rho(0,x_2;x_3,0)^2
q_{2,1}^{(4\leftrightarrow 1)}(0,x_2,x_3)
\,\mathrm{d}x_2\,\mathrm{d}x_3,
\\
\sigma_{1,(4\leftrightarrow 2),V_2}^{2}
&=
\int_{(\mathbf{R}^2)^2}
\rho(0,x_2;x_3,x_2)^2
q_{2,1}^{(4\leftrightarrow 2)}(0,x_2,x_3)
\,\mathrm{d}x_2\,\mathrm{d}x_3.
\end{align*}
Equivalently, these four quantities may be written compactly as
\[
\left(
\sigma_{1,(j\leftrightarrow i),V_2}^{2}
\right)_{(j,i)\in\mathcal{P}_2}.
\]

Proposition~\ref{prop:variance} below shows that the conditional variance
converges almost surely to
\[
\sigma_{V_2}^{2}
=
2
+
\frac{2}{3}
\left(
\sigma_{0,V_2}^{2}
+
\sum_{(j,i)\in\mathcal{P}_2}
\sigma_{1,(j\leftrightarrow i),V_2}^{2}
\right).
\]
The term \(2\) is the diagonal contribution, corresponding to
\((x_1,x_2)=(x_3,x_4)\). The remaining terms correspond respectively to
pairs of disjoint edges and pairs of distinct edges sharing one endpoint.

\begin{lemma}
\label{Le:finitevariances}
With the above notation, \(\sigma_{V_2}^{2}\) is finite.
\end{lemma}

\begin{prooft}{Lemma \ref{Le:finitevariances}}

We first prove that \(\sigma_{0,V_2}^{2}<\infty\). Let us write again
\[
\rho(x_1,x_2;x_3,x_4)
=
\mathrm{corr}
\left(
U_{x_1,x_2}^{(W)},U_{x_3,x_4}^{(W)}
\right).
\]
By symmetry between the two edges, it is enough to control the part of the
integral corresponding to
\[
\|x_4-x_3\|\leq \|x_2\|.
\]
The complementary part is treated in the same way after exchanging the two
edges and using translation invariance.

Set
\[
\ell_1=\|x_2\|,
\qquad
\ell_2=\|x_4-x_3\|.
\]
We therefore consider
\[
g(x_2,x_3,x_4)
=
\rho(0,x_2;x_3,x_4)^2
p_{2,1}(0,x_2,x_3,x_4)
\mathbb{I}[\ell_2\leq \ell_1].
\]
We show that \(g\) is integrable on \((\mathbf{R}^2)^3\).

Let \(\varepsilon\in(0,1/2)\), and let \(d_0\) be as in
Lemma~\ref{Le:bound:corr}\textit{(ii)}. We split the domain into three
regions.

\medskip
\noindent
\textbf{Case 1.}
Assume that
\[
\ell_1\leq \|x_3\|^{\varepsilon},
\qquad
\|x_3\|\geq d_0.
\]
By Lemma~\ref{Le:bound:corr}\textit{(ii)},
\[
|\rho(0,x_2;x_3,x_4)|
\leq
c\ell_1^{2-2H}\|x_3\|^{-2}.
\]
Moreover, by Lemma~\ref{Le:estimatepN}, with \(N=1\),
\[
p_{2,1}(0,x_2,x_3,x_4)
\leq
c(1+\ell_1^2)e^{-\frac{\pi}{4}\ell_1^2},
\]
on the region \(\ell_2\leq \ell_1\). Hence
\[
g(x_2,x_3,x_4)
\leq
c\,
\ell_1^{4-4H}
\|x_3\|^{4H-4}
(1+\ell_1^2)
e^{-\frac{\pi}{4}\ell_1^2}
\mathbb{I}[\ell_2\leq \ell_1].
\]
Integrating first with respect to \(x_4\), we use
\[
\int_{\mathbf{R}^2}
\mathbb{I}[\|x_4-x_3\|\leq \ell_1]\,\mathrm{d}x_4
=
\pi \ell_1^2.
\]
Thus the corresponding integral is bounded by
\[
c
\left(
\int_{\mathbf{R}^2}
(\ell_1^{6-4H}+\ell_1^{8-4H})
e^{-\frac{\pi}{4}\ell_1^2}
\,\mathrm{d}x_2
\right)
\left(
\int_{\mathbf{R}^2}
\|x_3\|^{4H-4}
\mathbb{I}[\|x_3\|\geq d_0]
\,\mathrm{d}x_3
\right).
\]
The first integral is finite because of the exponential factor. The second
one is finite if and only if
\[
4H-4<-2,
\]
which is equivalent to \(H<1/2\). Therefore the contribution of Case 1 is
finite.

\medskip
\noindent
\textbf{Case 2.}
Assume that \(\|x_3\|<d_0\). Since \(|\rho|\leq1\), Lemma
\ref{Le:estimatepN} gives
\[
g(x_2,x_3,x_4)
\leq
c(1+\ell_1^2)e^{-\frac{\pi}{4}\ell_1^2}
\mathbb{I}[\ell_2\leq \ell_1]
\mathbb{I}[\|x_3\|<d_0].
\]
After integration with respect to \(x_4\), this is bounded by
\[
c
\ell_1^2(1+\ell_1^2)e^{-\frac{\pi}{4}\ell_1^2}
\mathbb{I}[\|x_3\|<d_0],
\]
which is integrable with respect to \((x_2,x_3)\). Hence the contribution of
Case 2 is finite.

\medskip
\noindent
\textbf{Case 3.}
Assume that
\[
\ell_1>\|x_3\|^{\varepsilon}.
\]
Then
\[
\|x_3\|<\ell_1^{1/\varepsilon}.
\]
Again using \(|\rho|\leq1\) and Lemma~\ref{Le:estimatepN},
\[
g(x_2,x_3,x_4)
\leq
c(1+\ell_1^2)e^{-\frac{\pi}{4}\ell_1^2}
\mathbb{I}[\ell_2\leq \ell_1]
\mathbb{I}[\|x_3\|<\ell_1^{1/\varepsilon}].
\]
Integrating with respect to \(x_4\) and \(x_3\), the corresponding
contribution is bounded by
\[
c
\int_{\mathbf{R}^2}
\ell_1^2
(1+\ell_1^2)
\ell_1^{2/\varepsilon}
e^{-\frac{\pi}{4}\ell_1^2}
\,\mathrm{d}x_2,
\]
which is finite because of the exponential factor. This proves that
\(\sigma_{0,V_2}^{2}<\infty\).

It remains to consider the terms
\(\sigma_{1,(j\leftrightarrow i),V_2}^{2}\), for
\((j,i)\in\mathcal{P}_2\). We only give the argument once, since the four
cases are identical up to a relabelling of the two free edge vectors.
For such a term, after setting \(x_1=0\), the integral is over two free
vectors corresponding to two Delaunay edges sharing one endpoint. Denote
their lengths by \(\ell_a\) and \(\ell_b\), and set
\[
R=\max(\ell_a,\ell_b).
\]
By Lemma~\ref{Le:estimatepN}, with \(N=1\),
\[
q_{2,1}^{(j\leftrightarrow i)}
\leq
c(1+R^2)e^{-\frac{\pi}{4}R^2}.
\]
Since correlations are bounded by \(1\), the integrand is bounded by
\[
c(1+R^2)e^{-\frac{\pi}{4}R^2}.
\]

After the linear change of variables from the two free points to the two edge
vectors, say \(h_a,h_b\in\mathbb{R}^2\), the integral is bounded by
\[
c\int_{((\mathbf{R}^2)^2}
\left(1+\max\{\|h_a\|,\|h_b\|\}^2\right)
\exp\left\{-\frac{\pi}{4}\max\{\|h_a\|,\|h_b\|\}^2\right\}
\,dh_a\,dh_b .
\]
This is a four-dimensional integral with respect to the two vector variables
\(h_a\) and \(h_b\), and it is finite because the exponential decay in
\(\max\{\|h_a\|,\|h_b\|\}\) dominates any polynomial factor.
Therefore
\[
\sigma_{1,(j\leftrightarrow i),V_2}^{2}<\infty,
\qquad
(j,i)\in\mathcal{P}_2.
\]

We conclude that all the terms entering the definition of
\(\sigma_{V_2}^{2}\) are finite, and hence
\[
\sigma_{V_2}^{2}<\infty.
\]
\end{prooft}

We are now well equipped to prove  the following proposition.

\begin{proposition}
\label{prop:variance}
For \(H\in(0,1/2)\), the conditional second moment
\[
\mathbb{E}
\left[
\left(V_{2,N}^{(W)\prime}\right)^2
\,\middle|\,P_1
\right]
\]
converges, for almost every realization of \(P_1\), to
\[
\sigma_{V_2}^{2}
=
2+
\frac{2}{3}
\left(
\sigma_{0,V_2}^{2}
+
\sum_{(j,i)\in\mathcal P_2}
\sigma_{1,(j\leftrightarrow i),V_2}^{2}
\right).
\]
\end{proposition}

\begin{prooft}{Proposition \ref{prop:variance}}

We write, for shortness,
\[
\rho(x_1,x_2;x_3,x_4)
=
\mathrm{corr}
\left(
U_{x_1,x_2}^{(W)},U_{x_3,x_4}^{(W)}
\right).
\]
Recall that
\[
\mathbb{E}
\left[
\left(V_{2,N}^{(W)\prime}\right)^2
\,\middle|\,P_1
\right]
=
\frac{2}{|E_N'|}
\sum_{(x_1,x_2)\in E_N'}
\sum_{(x_3,x_4)\in E_N'}
\rho(x_1,x_2;x_3,x_4)^2.
\]
The diagonal contribution is equal to \(2\). We therefore have to identify
the limit of the non-diagonal contribution.

We introduce marks anchored at the first vertex of the first edge. For
\(x\in P_1\), define
\[
m_{0,x}
=
2
\sum_{\substack{x\sim x_2\\ x\preceq x_2}}
\sum_{\substack{x_3\sim x_4\\ x_3\preceq x_4}}
\rho(x,x_2;x_3,x_4)^2
\mathbb{I}
\left[
\{x,x_2\}\cap\{x_3,x_4\}=\emptyset
\right].
\]
This mark corresponds to the contribution of pairs of disjoint edges.

For pairs of distinct edges sharing one endpoint, we define four marks:
\begin{align*}
m_{1,(3\leftrightarrow 1),x}
&=
2
\sum_{\substack{x\sim x_2,\ x\sim x_4\\
x\preceq x_2,\ x\preceq x_4}}
\rho(x,x_2;x,x_4)^2
\mathbb{I}[x_2\neq x_4],
\\
m_{1,(3\leftrightarrow 2),x}
&=
2
\sum_{\substack{x\sim x_2,\ x_2\sim x_4\\
x\preceq x_2,\ x_2\preceq x_4}}
\rho(x,x_2;x_2,x_4)^2
\mathbb{I}[x\neq x_4],
\\
m_{1,(4\leftrightarrow 1),x}
&=
2
\sum_{\substack{x\sim x_2,\ x_3\sim x\\
x\preceq x_2,\ x_3\preceq x}}
\rho(x,x_2;x_3,x)^2
\mathbb{I}[x_2\neq x_3],
\\
m_{1,(4\leftrightarrow 2),x}
&=
2
\sum_{\substack{x\sim x_2,\ x_3\sim x_2\\
x\preceq x_2,\ x_3\preceq x_2}}
\rho(x,x_2;x_3,x_2)^2
\mathbb{I}[x\neq x_3].
\end{align*}
Each of these marks is non-negative and translation-covariant as a function
of the Poisson--Delaunay tessellation. By Lemma~\ref{Le:finitevariances},
their Palm expectations are finite. Since the homogeneous Poisson point
process is stationary and ergodic, the corresponding marked point processes
are ergodic. Hence the ergodic theorem gives, almost surely,
\[
\frac1N\sum_{x\in P_1\cap \mathbf C_N}m_{0,x}
\longrightarrow
\mathbb{E}
\left[
\sum_{x\in P_1\cap \mathbf C_1}m_{0,x}
\right],
\]
and similarly for the four marks \(m_{1,(j\leftrightarrow i),x}\).

By the Slivnyak--Mecke formula and translation invariance,
\[
\mathbb{E}
\left[
\sum_{x\in P_1\cap \mathbf C_1}m_{0,x}
\right]
=
2\sigma_{0,V_2}^{2}.
\]
Likewise,
\[
\mathbb{E}
\left[
\sum_{x\in P_1\cap \mathbf C_1}
m_{1,(j\leftrightarrow i),x}
\right]
=
2\sigma_{1,(j\leftrightarrow i),V_2}^{2},
\qquad
(j,i)\in\mathcal P_2.
\]
Consequently,
\begin{equation}
\frac{2}{N}
\sum_{\substack{x_1\in\mathbf C_N\\
x_1\sim x_2,\ x_1\preceq x_2}}
\sum_{\substack{x_3\sim x_4\\
x_3\preceq x_4}}
\rho(x_1,x_2;x_3,x_4)^2
\mathbb{I}
\left[
(x_1,x_2)\neq(x_3,x_4)
\right]
\longrightarrow
2
\left(
\sigma_{0,V_2}^{2}
+
\sum_{(j,i)\in\mathcal P_2}
\sigma_{1,(j\leftrightarrow i),V_2}^{2}
\right)
\label{eq:non-diagonal-all-space}
\end{equation}
almost surely.

It remains to justify that replacing the second sum over all Delaunay edges
by a sum over edges whose first endpoint belongs to \(\mathbf C_N\) does not
change the limit. That is, we have to prove that
\[
B_N
:=
\frac{2}{N}
\sum_{\substack{x_1\in\mathbf C_N\\
x_1\sim x_2,\ x_1\preceq x_2}}
\sum_{\substack{x_3\notin\mathbf C_N\\
x_3\sim x_4,\ x_3\preceq x_4}}
\rho(x_1,x_2;x_3,x_4)^2
\longrightarrow 0
\]
almost surely.

We use a standard truncation argument. For \(d>0\), let
\[
\mathbf A_N(d)=(-N^{1/2}/2+d,N^{1/2}/2-d]^2
\]
and let \(D_d(x)\) be the square centered at \(x\), with side length \(2d\)
and sides parallel to the coordinate axes. Define the truncated contribution
\[
S_N(d)
=
\frac{2}{N}
\sum_{\substack{x_1\in\mathbf A_N(d)\\
x_1\sim x_2,\ x_1\preceq x_2}}
\sum_{\substack{x_3\in D_d(x_1)\\
x_3\sim x_4,\ x_3\preceq x_4}}
\rho(x_1,x_2;x_3,x_4)^2.
\]
For fixed \(d\), the same ergodic argument as above yields
\[
S_N(d)
\longrightarrow
6+
2
\left(
\sigma_{0,V_2}^{2}(d)
+
\sum_{(j,i)\in\mathcal P_2}
\sigma_{1,(j\leftrightarrow i),V_2}^{2}(d)
\right)
=: L(d),
\]
where the truncated constants are obtained by adding the constraint
\(x_3\in D_d(0)\) in the disjoint-edge contribution and the analogous
constraint in the one-common-endpoint contributions.

By monotone convergence and Lemma~\ref{Le:finitevariances},
\[
L(d)
\longrightarrow
6+
2
\left(
\sigma_{0,V_2}^{2}
+
\sum_{(j,i)\in\mathcal P_2}
\sigma_{1,(j\leftrightarrow i),V_2}^{2}
\right)
=:L
\]
as \(d\to\infty\).

On the other hand, for fixed \(d\) and \(N\) large enough, if
\(x_1\in\mathbf A_N(d)\) and \(x_3\in D_d(x_1)\), then
\(x_3\in\mathbf C_N\). Therefore
\[
S_N(d)
\leq
\frac{2}{N}
\sum_{\substack{x_1\in\mathbf C_N\\
x_1\sim x_2,\ x_1\preceq x_2}}
\sum_{\substack{x_3\in\mathbf C_N\\
x_3\sim x_4,\ x_3\preceq x_4}}
\rho(x_1,x_2;x_3,x_4)^2
\leq
S_N(\infty),
\]
where \(S_N(\infty)\) denotes the same quantity with the second sum taken
over all Delaunay edges. We already know from \eqref{eq:non-diagonal-all-space}
and the diagonal contribution that
\[
S_N(\infty)\longrightarrow L
\]
almost surely. The preceding inequalities and the convergence
\(L(d)\to L\) imply
\[
\frac{2}{N}
\sum_{\substack{x_1\in\mathbf C_N\\
x_1\sim x_2,\ x_1\preceq x_2}}
\sum_{\substack{x_3\in\mathbf C_N\\
x_3\sim x_4,\ x_3\preceq x_4}}
\rho(x_1,x_2;x_3,x_4)^2
\longrightarrow
L.
\]
Hence the contribution with \(x_3\notin\mathbf C_N\) is negligible:
\[
B_N\longrightarrow 0.
\]

Combining this boundary estimate with the diagonal contribution gives
\[
\frac{2}{N}
\sum_{(x_1,x_2)\in E_N'}
\sum_{(x_3,x_4)\in E_N'}
\rho(x_1,x_2;x_3,x_4)^2
\longrightarrow
6+
2
\left(
\sigma_{0,V_2}^{2}
+
\sum_{(j,i)\in\mathcal P_2}
\sigma_{1,(j\leftrightarrow i),V_2}^{2}
\right).
\]
Finally, since
\[
\frac{|E_N'|}{N}\longrightarrow 3
\]
almost surely, we obtain
\begin{align*}
\mathbb{E}
\left[
\left(V_{2,N}^{(W)\prime}\right)^2
\,\middle|\,P_1
\right]
&=
\frac{2}{|E_N'|}
\sum_{(x_1,x_2)\in E_N'}
\sum_{(x_3,x_4)\in E_N'}
\rho(x_1,x_2;x_3,x_4)^2
\\
&\longrightarrow
\frac{1}{3}
\left[
6+
2
\left(
\sigma_{0,V_2}^{2}
+
\sum_{(j,i)\in\mathcal P_2}
\sigma_{1,(j\leftrightarrow i),V_2}^{2}
\right)
\right]
\\
&=
2+
\frac{2}{3}
\left(
\sigma_{0,V_2}^{2}
+
\sum_{(j,i)\in\mathcal P_2}
\sigma_{1,(j\leftrightarrow i),V_2}^{2}
\right).
\end{align*}
This proves the proposition.
\end{prooft}

\subsection{Vanishing of the first contraction}
\label{Sect_Cond_c}

We now prove that the first contraction of \(f_{N,2}\) vanishes. Since
\(V_{2,N}^{(W)\prime}=I_{2}(f_{N,2})\), the only contraction to control is
the contraction of order one. Recall that
\[
f_{N,2}
=
\frac{1}{\sqrt{|e_N'|}}
\sum_{k\in e_N'}
\varepsilon_k^{\otimes 2}.
\]
Thus
\[
f_{N,2}\otimes_1 f_{N,2}
=
\frac{1}{|e_N'|}
\sum_{k,l\in e_N'}
\langle \varepsilon_k,\varepsilon_l\rangle_{\mathfrak H}
\varepsilon_k\otimes \varepsilon_l.
\]
Consequently,
\begin{equation}
\left\|
f_{N,2}\otimes_1 f_{N,2}
\right\|_{\mathfrak H^{\otimes 2}}^{2}
=
\frac{1}{|e_N'|^{2}}
\sum_{k,l,i,j\in e_N'}
\rho_{k,l}\rho_{i,j}\rho_{k,i}\rho_{l,j},
\label{eq:contraction-expansion}
\end{equation}
where
\[
\rho_{a,b}
=
\mathrm{corr}(U^{(a)},U^{(b)}).
\]
It is enough to prove that the absolute value of the right-hand side of
\eqref{eq:contraction-expansion} converges to \(0\) almost surely.

Let \(d_{k,l}\) denote the distance between the first point of the edge
\(\varphi(k)\) and the first point of the edge \(\varphi(l)\). Also, let
\(\|e_m\|\) denote the length of the edge \(\varphi(m)\). Let
\(\varepsilon\in(0,1/2)\), and let \(d_0\) be as in
Lemma~\ref{Le:bound:corr}\textit{(ii)}.

We first consider the most delicate case, namely when
\[
d_{k,l},d_{i,j},d_{k,i},d_{l,j}
\geq
\max\{d_0,\|e_i\|^{1/\varepsilon}\},
\]
and where, without loss of generality,
\[
\|e_i\|\geq \|e_j\|,\|e_k\|,\|e_l\|.
\]
The remaining configurations, where at least one of these distances is not
large compared with the edge lengths, are handled in the same way by using
the trivial bound \(|\rho_{a,b}|\leq 1\) for the corresponding correlation
and by restricting one of the indices to a local neighborhood. These terms
are of smaller order than the one treated below.

By Lemma~\ref{Le:bound:corr}\textit{(ii)}, it is enough to prove that
\begin{multline}
\frac{1}{|e_N'|^{2}}
\sum_{k,l,i,j\in e_N'}
d_{k,l}^{2H-2}\mathbb{I}[d_{k,l}\geq d_0]\,
d_{i,j}^{2H-2}\mathbb{I}[d_{i,j}\geq d_0]
\\
\times
d_{k,i}^{2H-2}\mathbb{I}[d_{k,i}\geq d_0]\,
d_{l,j}^{2H-2}\mathbb{I}[d_{l,j}\geq d_0]\,
\|e_i\|^{8(2-2H)}
\longrightarrow 0
\label{eq:contraction-bound-to-prove}
\end{multline}
almost surely.

We split the proof according to whether the longest edge involved in the
four correlations is larger or smaller than \(N^{\varepsilon_0}\), where
\(\varepsilon_0>0\) will be chosen sufficiently small.

\medskip
\noindent
\textbf{Case 1: large edges.}
Assume first that \(\|e_i\|\geq N^{\varepsilon_0}\). Set
\[
p=8(2-2H).
\]
The contribution of such terms to the left-hand side of
\eqref{eq:contraction-bound-to-prove} is bounded by
\[
c|e_N'|
\sum_{i\in e_N'}
\|e_i\|^{p}
\mathbb{I}[\|e_i\|\geq N^{\varepsilon_0}],
\]
for a constant \(c>0\).

We now prove that this upper bound converges to \(0\) almost surely. The
proof of Lemma~\ref{Le:estimatepN} gives the following single-edge estimate:
there exists \(c>0\) such that, for any two distinct points \(x,z\),
\[
\mathbb{P}
\left[
x\sim z
\text{ in }
\operatorname{Del}(P_1\cup\{x,z\})
\right]
\leq
c(1+\|z-x\|^2)
\exp\left\{
-\frac{\pi}{4}\|z-x\|^2
\right\}.
\]
Therefore, by the Slivnyak--Mecke formula,
\begin{align*}
\mathbb{E}
\left[
\sum_{i\in e_N'}
\|e_i\|^p
\mathbb{I}[\|e_i\|\geq N^{\varepsilon_0}]
\right]
&\leq
c
\int_{\mathbf C_N}
\int_{\mathbf R^2}
\|z-x\|^p
(1+\|z-x\|^2)
e^{-\frac{\pi}{4}\|z-x\|^2}
\\
&\hspace{3cm}\times
\mathbb{I}[\|z-x\|\geq N^{\varepsilon_0}]
\,\mathrm{d}z\,\mathrm{d}x
\\
&\leq
cN
\int_{N^{\varepsilon_0}}^\infty
r^{p+1}(1+r^2)e^{-\frac{\pi}{4}r^2}
\,\mathrm{d}r
\\
&\leq
cN
\exp\{-cN^{2\varepsilon_0}\},
\end{align*}
for a possibly different constant \(c>0\). Hence
\[
\mathbb{E}
\left[
N
\sum_{i\in e_N'}
\|e_i\|^p
\mathbb{I}[\|e_i\|\geq N^{\varepsilon_0}]
\right]
\leq
cN^2\exp\{-cN^{2\varepsilon_0}\}.
\]
The right-hand side is summable in \(N\). By Markov's inequality and the
Borel--Cantelli lemma,
\[
N
\sum_{i\in e_N'}
\|e_i\|^p
\mathbb{I}[\|e_i\|\geq N^{\varepsilon_0}]
\longrightarrow 0
\qquad\text{a.s.}
\]
Since \(|e_N'|/N\to 3\) almost surely, we also have
\[
|e_N'|
\sum_{i\in e_N'}
\|e_i\|^p
\mathbb{I}[\|e_i\|\geq N^{\varepsilon_0}]
\longrightarrow 0
\qquad\text{a.s.}
\]
This proves that the contribution of large edges is negligible.

\medskip
\noindent
\textbf{Case 2: moderate edges.}
Assume now that
\[
\|e_i\|\leq N^{\varepsilon_0}.
\]
Set
\[
\widetilde{\varepsilon}_0=8\varepsilon_0(2-2H).
\]
Then
\[
\|e_i\|^{8(2-2H)}
\leq
N^{\widetilde{\varepsilon}_0}.
\]
It remains to show that
\begin{equation}
\frac{N^{\widetilde{\varepsilon}_0}}{|e_N'|^{2}}
\sum_{k,l,i,j\in e_N'}
d_{k,l}^{2H-2}\mathbb{I}[d_{k,l}\geq d_0]\,
d_{i,j}^{2H-2}\mathbb{I}[d_{i,j}\geq d_0]
\times
d_{k,i}^{2H-2}\mathbb{I}[d_{k,i}\geq d_0]\,
d_{l,j}^{2H-2}\mathbb{I}[d_{l,j}\geq d_0]
\longrightarrow 0
\label{eq:moderate-edge-contribution}
\end{equation}
almost surely.

Since \(2H-2<0\), we have
\begin{equation*}
d_{k,l}^{2H-2}\mathbb{I}[d_{k,l}\geq d_0]\,
d_{k,i}^{2H-2}\mathbb{I}[d_{k,i}\geq d_0]
\leq
d_{k,l}^{2(2H-2)}\mathbb{I}[d_{k,l}\geq d_0]
+
d_{k,i}^{2(2H-2)}\mathbb{I}[d_{k,i}\geq d_0].
\end{equation*}
By symmetry, it is enough to control
\[
\sum_{k,l,i,j\in e_N'}
d_{k,l}^{2(2H-2)}\mathbb{I}[d_{k,l}\geq d_0]\,
d_{i,j}^{2H-2}\mathbb{I}[d_{i,j}\geq d_0]\,
d_{l,j}^{2H-2}\mathbb{I}[d_{l,j}\geq d_0].
\]

We first prove that, for every \(\eta>0\),
\[
N^{-\eta}
\sup_{l\in e_N'}
\sum_{k}
d_{k,l}^{2(2H-2)}
\mathbb{I}[d_{k,l}\geq d_0]
\longrightarrow 0
\qquad\text{a.s.}
\]
It is enough to show that
\[
N^{-\eta}
\sup_{y\in P_1\cap\mathbf C_N}
\sum_{x\in P_1}
\|y-x\|^{2(2H-2)}
\mathbb{I}[\|y-x\|\geq d_0]
\longrightarrow 0
\qquad\text{a.s.}
\]
Let \(\delta>0\). By the union bound and the Slivnyak--Mecke formula,
\begin{align*}
&\mathbb{P}
\left[
N^{-\eta}
\sup_{y\in P_1\cap\mathbf C_N}
\sum_{x\in P_1}
\|y-x\|^{2(2H-2)}
\mathbb{I}[\|y-x\|\geq d_0]
>
\delta
\right]
\\
&\leq
\mathbb{E}
\left[
\sum_{y\in P_1\cap\mathbf C_N}
\mathbb{I}
\left[
N^{-\eta}
\sum_{x\in P_1}
\|y-x\|^{2(2H-2)}
\mathbb{I}[\|y-x\|\geq d_0]
>
\delta
\right]
\right]
\\
&=
N
\mathbb{P}
\left[
N^{-\eta}
\sum_{x\in P_1}
\|x\|^{2(2H-2)}
\mathbb{I}[\|x\|\geq d_0]
>
\delta
\right].
\end{align*}
By Chernoff's inequality,
\begin{align*}
&\mathbb{P}
\left[
N^{-\eta}
\sum_{x\in P_1\cap B(0,\sqrt{2}%
N^{1/2})}
\|x\|^{2(2H-2)}
\mathbb{I}[\|x\|\geq d_0]
>
\delta
\right]
\\
&\leq
\exp(-\delta N^\eta)
\mathbb{E}
\left[
\exp
\left(
\sum_{x\in P_1}
\|x\|^{2(2H-2)}
\mathbb{I}[\|x\|\geq d_0]
\right)
\right].
\end{align*}
By the exponential formula for Poisson point processes,
\begin{equation*}
\mathbb{E}
\left[
\exp
\left(
\sum_{x\in P_1}
\|x\|^{2(2H-2)}
\mathbb{I}[\|x\|\geq d_0]
\right)
\right]
=
\exp
\left(
\int_{\mathbf R^2}
\left[
\exp
\left(
\|x\|^{2(2H-2)}
\mathbb{I}[\|x\|\geq d_0]
\right)
-1
\right]
\,\mathrm{d}x
\right).
\end{equation*}
The integral is finite because, as \(\|x\|\to\infty\),
\[
\exp(\|x\|^{2(2H-2)})-1
\sim
\|x\|^{2(2H-2)},
\]
and
\[
\int_{\|x\|\geq d_0}
\|x\|^{2(2H-2)}
\,\mathrm{d}x
<\infty
\]
when \(2H<1\). Therefore
\[
\mathbb{P}
\left[
N^{-\eta}
\sup_{y\in P_1\cap\mathbf C_N}
\sum_{x\in P_1}
\|y-x\|^{2(2H-2)}
\mathbb{I}[\|y-x\|\geq d_0]
>
\delta
\right]
\]
is summable in \(N\). The Borel--Cantelli lemma yields
\[
N^{-\eta}
\sup_{l\in e_N'}
\sum_{k}
d_{k,l}^{2(2H-2)}
\mathbb{I}[d_{k,l}\geq d_0]
\longrightarrow 0
\qquad\text{a.s.}
\]

Thus, for every \(\eta_0>0\), almost surely and for \(N\) large enough,
\[
\sup_{l\in e_N'}
\sum_{k}
d_{k,l}^{2(2H-2)}
\mathbb{I}[d_{k,l}\geq d_0]
\leq
N^{\eta_0/2}.
\]
It follows that
\begin{multline*}
\sum_{k,l,i,j\in e_N'}
d_{k,l}^{2(2H-2)}\mathbb{I}[d_{k,l}\geq d_0]\,
d_{i,j}^{2H-2}\mathbb{I}[d_{i,j}\geq d_0]\,
d_{l,j}^{2H-2}\mathbb{I}[d_{l,j}\geq d_0]
\\
\leq
N^{\eta_0/2}
\sum_{l,i,j\in e_N'}
d_{i,j}^{2H-2}\mathbb{I}[d_{i,j}\geq d_0]\,
d_{l,j}^{2H-2}\mathbb{I}[d_{l,j}\geq d_0]
\qquad\text{a.s.}
\end{multline*}
Moreover,
\[
\sum_{l,i,j\in e_N'}
d_{i,j}^{2H-2}\mathbb{I}[d_{i,j}\geq d_0]\,
d_{l,j}^{2H-2}\mathbb{I}[d_{l,j}\geq d_0]
\leq
\sum_{j\in e_N'}
\left(
\sum_{i\in B(j,\sqrt{2}N^{1/2})}
d_{i,j}^{2H-2}\mathbb{I}[d_{i,j}\geq d_0]
\right)^2.
\]
Here \(B(j,\sqrt{2}N^{1/2})\) denotes the set of indices \(i\in e_N'\)
whose first endpoint is at distance at most \(\sqrt{2}N^{1/2}\) from the
first endpoint of the edge \(\varphi(j)\).

Choose \(\varepsilon_0>0\) and \(\eta_0>0\) such that
\[
\widetilde{\varepsilon}_0+\frac{\eta_0}{2}<1-2H.
\]
It remains to prove that, for some \(\eta>0\) satisfying
\[
\widetilde{\varepsilon}_0+\frac{\eta_0}{2}+2\eta<1-2H,
\]
one has
\[
N^{-(2H/2+\eta)}
\sup_{y\in P_1\cap\mathbf C_N}
\sum_{x\in P_1\cap B(y,\sqrt{2}N^{1/2})}
\|y-x\|^{2H-2}
\mathbb{I}[\|y-x\|\geq d_0]
\longrightarrow 0
\qquad\text{a.s.}
\]
Indeed, since \(|P_1\cap\mathbf C_N|/N\to 1\) a.s., this bound implies
\begin{equation*}
\frac{N^{\widetilde{\varepsilon}_0+\eta_0/2}}{|e_N'|^2}
\sum_{j\in e_N'}
\left(
\sum_{i\in B(j,\sqrt{2}N^{1/2})}
d_{i,j}^{2H-2}\mathbb{I}[d_{i,j}\geq d_0]
\right)^2
=
O\left(
N^{\widetilde{\varepsilon}_0+\eta_0/2+2H+2\eta-1}
\right)
\longrightarrow 0
\qquad\text{a.s.}
\end{equation*}

It remains to establish the displayed maximal bound. Let \(\delta>0\). By
the union bound and the Slivnyak--Mecke formula,
\begin{align*}
&\mathbb{P}
\left[
N^{-(H+\eta)}
\sup_{y\in P_1\cap\mathbf C_N}
\sum_{x\in P_1\cap B(y,\sqrt{2}N^{1/2})}
\|y-x\|^{2H-2}
\mathbb{I}[\|y-x\|\geq d_0]
>
\delta
\right]
\\
&\leq
N
\mathbb{P}
\left[
N^{-(H+\eta)}
\sum_{x\in P_1\cap B(0,\sqrt{2}N^{1/2})}
\|x\|^{2H-2}
\mathbb{I}[\|x\|\geq d_0]
>
\delta
\right].
\end{align*}
By Chernoff's inequality,
\begin{align*}
&\mathbb{P}
\left[
N^{-(H+\eta)}
\sum_{x\in P_1\cap B(0,\sqrt{2}N^{1/2})}
\|x\|^{2H-2}
\mathbb{I}[\|x\|\geq d_0]
>
\delta
\right]
\\
&\leq
\exp(-\delta N^{H+\eta})
\mathbb{E}
\left[
\exp
\left(
\sum_{x\in P_1\cap B(0,\sqrt{2}N^{1/2})}
\|x\|^{2H-2}
\mathbb{I}[\|x\|\geq d_0]
\right)
\right].
\end{align*}
Let
\[
\mathcal N
=
\left|P_1\cap B(0,\sqrt{2}N^{1/2})\right|.
\]
Then \(\mathcal N\) is Poisson distributed with parameter \(2\pi N\).
Conditionally on \(\mathcal N\), the points are independent and uniformly
distributed in \(B(0,\sqrt{2}N^{1/2})\). If \(R\) denotes the distance of
such a uniform point to the origin, then
\[
\mathbb{P}(R\leq r)
=
\frac{r^2}{2N},
\qquad
0\leq r\leq \sqrt{2N}.
\]
Therefore Lemma~\ref{Le_Bound_Exp_R} applies with \(2N\) in place of \(N\),
and for large \(N\),
\[
\mathbb{E}
\left[
e^{R^{2H-2}\mathbb{I}[R\geq d_0]}
\right]
\leq
1+c_1(\sqrt{2N})^{2H-2}+c_2(2N)^{-1}.
\]
Using the moment generating function of the Poisson distribution, we obtain
\begin{align*}
&\mathbb{E}
\left[
\exp
\left(
\sum_{x\in P_1\cap B(0,\sqrt{2}N^{1/2})}
\|x\|^{2H-2}
\mathbb{I}[\|x\|\geq d_0]
\right)
\right]
\\
&\leq
\exp
\left\{
2\pi N
\left(
c_1(\sqrt{2N})^{2H-2}+c_2(2N)^{-1}
\right)
\right\}
\\
&\leq
\exp\left\{
C N^{H}+C
\right\},
\end{align*}
for some constant \(C>0\). Consequently,
\begin{align*}
&\mathbb{P}
\left[
N^{-(H+\eta)}
\sum_{x\in P_1\cap B(0,\sqrt{2}N^{1/2})}
\|x\|^{2H-2}
\mathbb{I}[\|x\|\geq d_0]
>
\delta
\right]
\\
&\leq
\exp\left\{
-\delta N^{H+\eta}
+
C N^{H}
+
C
\right\}.
\end{align*}
After multiplication by the prefactor \(N\) coming from the union bound,
the resulting sequence is still summable in \(N\), since \(\eta>0\). The
Borel--Cantelli lemma gives
\[
N^{-(H+\eta)}
\sup_{y\in P_1\cap\mathbf C_N}
\sum_{x\in P_1\cap B(y,\sqrt{2}N^{1/2})}
\|y-x\|^{2H-2}
\mathbb{I}[\|y-x\|\geq d_0]
\longrightarrow 0
\qquad\text{a.s.}
\]
This proves \eqref{eq:moderate-edge-contribution}.

Combining the estimates for large and moderate edges, we conclude that
\[
\left\|
f_{N,2}\otimes_1 f_{N,2}
\right\|_{\mathfrak H^{\otimes 2}}^{2}
\longrightarrow 0
\qquad\text{a.s.}
\]
This proves the contraction condition required for the second-chaos normal
approximation criterion.

\section{Simulation study}
\label{sec:simulations}

We illustrate the central limit theorem of Theorem~\ref{th:CLTgaussian}
through a Monte Carlo simulation study. The aim is to examine the
finite-sample behavior of the edge-based statistic \(V^{(W)}_{2,N}\) when
the fractional Brownian field is observed at random spatial locations and
increments are computed along the edges of the associated Delaunay
triangulation.

\begin{figure}[t]
\centering
\includegraphics[width=0.7\textwidth]{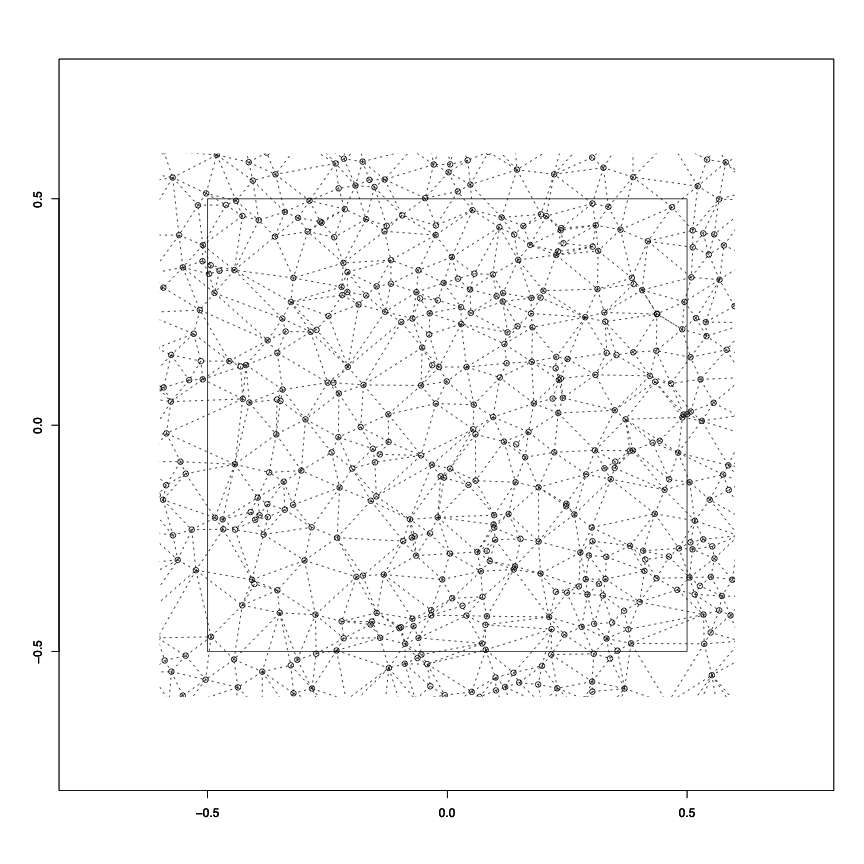}
\caption{A realization of a homogeneous Poisson point process in
\(\mathbf{C}=(-1/2,1/2]^{2}\) and the associated Delaunay triangulation.}
\label{fig:delaunay}
\end{figure}

We generate a homogeneous Poisson point process with intensity \(N=300\) in
the unit square $\mathbf{C}=(-1/2,1/2]^2$, so that the expected number of points in \(\mathbf{C}\) is equal to \(300\).
A typical realization of the point pattern, together with its Delaunay
triangulation, is displayed in Figure~\ref{fig:delaunay}. The figure
illustrates the irregular geometry induced by the random sampling scheme:
both the lengths and orientations of the Delaunay edges vary substantially,
in contrast with the regular-grid framework usually considered for
quadratic variations of Gaussian fields.

Conditionally on this realization of the Poisson point process and its
Delaunay triangulation, we simulate an isotropic fractional Brownian field
\(W\) with covariance function \eqref{eq:defcovariance}. The scale parameter
is set to \(\sigma=1\), without loss of generality. For each Delaunay edge
\((x_1,x_2)\), we compute the normalized increment
\[
U^{(W)}_{x_1,x_2}
 =
 \|x_2-x_1\|^{-H}
 \left(W(x_2)-W(x_1)\right),
\]
and form the statistic
\[
V^{(W)}_{2,N}
 =
 \frac{1}{\sqrt{|E_N|}}
 \sum_{(x_1,x_2)\in E_N}
 \left\{
   \left(U^{(W)}_{x_1,x_2}\right)^2 - 1
 \right\},
\]
where \(E_N\) denotes the set of Delaunay edges whose lexicographically first
endpoint belongs to \(\mathbf{C}\).

\begin{figure}[t]
\centering
\includegraphics[width=\textwidth]{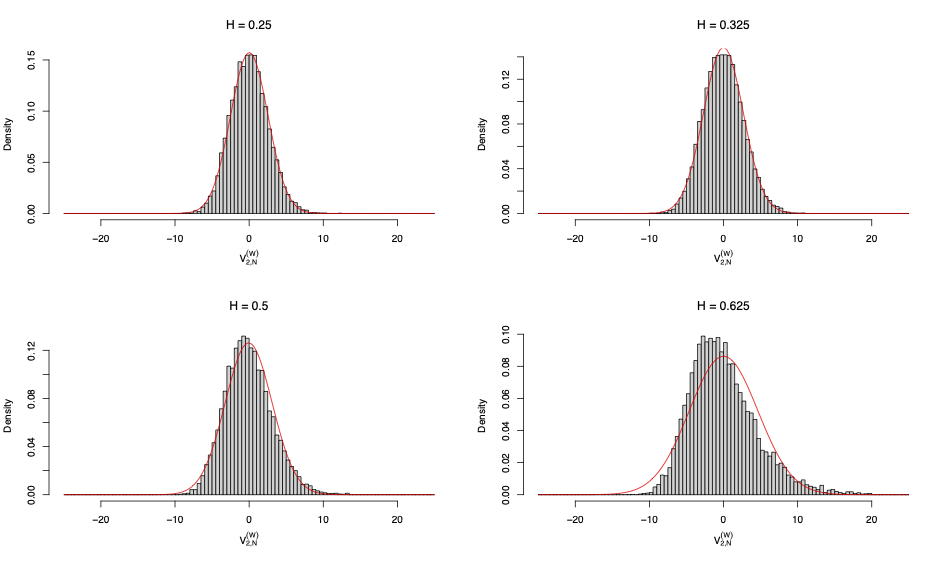}
\caption{Empirical distributions of the statistic \(V^{(W)}_{2,N}\), based
on \(10\,000\) independent simulations of the fractional Brownian field for
a fixed Poisson--Delaunay triangulation with approximately \(300\) points in
\(\mathbf{C}=(-1/2,1/2]^2\), and for several values of
\(H\). The red curves correspond to Gaussian densities with matching
empirical mean and variance.}
\label{fig:histV2}
\end{figure}

For each value of the parameter \(H\), we perform \(10\,000\)
independent simulations of the fractional Brownian field, keeping the
underlying Poisson point pattern fixed. Figure~\ref{fig:histV2} displays the
empirical histograms of \(V^{(W)}_{2,N}\) for several values of \(2H\).
The red curves correspond to Gaussian densities with empirical mean and
variance matched to the simulated samples.

For \(H<1/2\), the empirical distributions are close to the Gaussian
shape predicted by Theorem~\ref{th:CLTgaussian}, even for this moderate
number of points. By contrast, for the boundary case \(H=0.5\) and for
\(H=0.625\), which lies outside the range covered by the theorem, the
histograms exhibit visible departures from normality. In particular, the
empirical distributions become more skewed and display heavier tails. This
behavior is consistent with the theoretical restriction \(2H<1\), or
equivalently \(H<1/2\), under which the Gaussian limit is established.

Overall, the simulation study supports the relevance of the Gaussian
approximation in the admissible parameter range and illustrates the change
of asymptotic behavior that occurs beyond the critical threshold
\(H=0.5\).

\begin{appendices}

\section{Proof of Theorem \protect\ref{th:CLTgaussian} for \(V_{3,N}^{(W)}\)}
\label{sec:proofpairsedges}

We give a proof sketch, since the argument is the same as for the
edge-based statistic \(V_{2,N}^{(W)}\), up to heavier notation. The only
additional point is that each Delaunay triangle contributes to two
orthogonalized normalized increments instead of one normalized edge
increment.

By the same scaling argument as in Section~\ref{sec:proofedges},
\(V_{3,N}^{(W)}\) has the same distribution as
\[
V_{3,N}^{(W)\prime}
=
\frac{1}{\sqrt{|DT_N'|}}
\sum_{(x_1,x_2,x_3)\in DT_N'}
\left[
\left(
\begin{array}{cc}
U_{x_1,x_2}^{(W)} & U_{x_1,x_3}^{(W)}
\end{array}
\right)
\left(
\begin{array}{cc}
1 & R_{x_1,x_2,x_3}\\
R_{x_1,x_2,x_3} & 1
\end{array}
\right)^{-1}
\left(
\begin{array}{c}
U_{x_1,x_2}^{(W)}\\
U_{x_1,x_3}^{(W)}
\end{array}
\right)
-2
\right],
\]
where \(DT_N'=DT_{1,\mathbf C_N}\), with
\[
\mathbf C_N=(-N^{1/2}/2,N^{1/2}/2]^2.
\]
Thus \(DT_N'\) is the set of ordered Delaunay triangles
\((x_1,x_2,x_3)\) of \(\operatorname{Del}(P_1)\) such that
\[
\Delta(x_1,x_2,x_3)\in \operatorname{Del}(P_1),
\qquad
x_1\in\mathbf C_N,
\qquad
x_1\preceq x_2\preceq x_3.
\]
Moreover,
\[
\frac{|DT_N'|}{N}
\longrightarrow 2
\qquad\text{a.s.},
\]
because the cell intensity of the Poisson--Delaunay triangulation is equal
to \(2\).

For a set of three points
\[
\Delta:=\Delta(x_{1},x_{2},x_{3})=(x_1,x_2,x_3)\in DT_N',
\]
write
\[
R_\Delta=R_{x_1,x_2,x_3}.
\]
Define
\[
Z_{\Delta,1}
=
(1-R_\Delta^2)^{-1/2}
\left(
U_{x_1,x_2}^{(W)}
-
R_\Delta U_{x_1,x_3}^{(W)}
\right),
\qquad
Z_{\Delta,2}
=
U_{x_1,x_3}^{(W)}.
\]
Then \(Z_{\Delta,1}\) and \(Z_{\Delta,2}\) are standard Gaussian random
variables and
\[
\mathrm{corr}(Z_{\Delta,1},Z_{\Delta,2})=0.
\]
Consequently,
\[
V_{3,N}^{(W)\prime}
=
\frac{1}{\sqrt{|DT_N'|}}
\sum_{\Delta\in DT_N'}
\left[
  (Z_{\Delta,1}^{2}-1)
  +(Z_{\Delta,2}^{2}-1)
\right].
\]

Conditionally on \(P_1\), the family
\[
\{Z_{\Delta,a}:\Delta\in DT_N',\ a=1,2\}
\]
is a centered Gaussian family with unit variances. Hence there exist a real
separable Hilbert space \(\mathfrak H\), an isonormal Gaussian process
\(\{X(h):h\in\mathfrak H\}\), and vectors
\(\eta_{\Delta,a}\in\mathfrak H\), \(a=1,2\), such that
\[
Z_{\Delta,a}=X(\eta_{\Delta,a}),
\qquad
\|\eta_{\Delta,a}\|_{\mathfrak H}=1,
\]
and
\[
\langle \eta_{\Delta,a},\eta_{\Delta',b}\rangle_{\mathfrak H}
=
\mathrm{corr}(Z_{\Delta,a},Z_{\Delta',b}).
\]
Setting
\[
f_{N,3}
=
\frac{1}{\sqrt{|DT_N'|}}
\sum_{\Delta\in DT_N'}
\left(
\eta_{\Delta,1}^{\otimes 2}
+
\eta_{\Delta,2}^{\otimes 2}
\right),
\]
we obtain
\[
V_{3,N}^{(W)\prime}=I_2(f_{N,3}).
\]
Thus the proof again reduces to the second-chaos normal approximation
criterion. We have to prove the convergence of the conditional variance and
the vanishing of the first contraction.

We now record the correlations that enter the proof. Let
\[
\Delta=(x_1,x_2,x_3),
\qquad
\Delta'=(x_4,x_5,x_6),
\]
and write
\[
R_\Delta=R_{x_1,x_2,x_3},
\qquad
R_{\Delta'}=R_{x_4,x_5,x_6}.
\]
For shortness, set
\[
\rho_{ab,cd}
=
\mathrm{corr}
\left(
U_{x_a,x_b}^{(W)},U_{x_c,x_d}^{(W)}
\right).
\]
Then
\begin{align*}
\mathrm{corr}(Z_{\Delta,1},Z_{\Delta',1})
&=
(1-R_\Delta^2)^{-1/2}
(1-R_{\Delta'}^2)^{-1/2}
\\
&\quad\times
\left[
\rho_{12,45}
-
R_{\Delta'}\rho_{12,46}
-
R_\Delta\rho_{13,45}
+
R_\Delta R_{\Delta'}\rho_{13,46}
\right],
\\[1ex]
\mathrm{corr}(Z_{\Delta,1},Z_{\Delta',2})
&=
(1-R_\Delta^2)^{-1/2}
\left[
\rho_{12,46}
-
R_\Delta\rho_{13,46}
\right],
\\[1ex]
\mathrm{corr}(Z_{\Delta,2},Z_{\Delta',1})
&=
(1-R_{\Delta'}^2)^{-1/2}
\left[
\rho_{13,45}
-
R_{\Delta'}\rho_{13,46}
\right],
\\[1ex]
\mathrm{corr}(Z_{\Delta,2},Z_{\Delta',2})
&=
\rho_{13,46}.
\end{align*}
In particular, all correlations between the variables \(Z_{\Delta,a}\) and
\(Z_{\Delta',b}\) are finite linear combinations of correlations between
normalized increments of the form studied in Section~\ref{sec:proofedges},
with coefficients depending only on the local geometries of the two
triangles.

The conditional variance is
\[
\mathbb{E}
\left[
\left(V_{3,N}^{(W)\prime}\right)^2
\,\middle|\,P_1
\right]
=
\frac{2}{|DT_N'|}
\sum_{\Delta,\Delta'\in DT_N'}
\sum_{a,b=1}^{2}
\mathrm{corr}(Z_{\Delta,a},Z_{\Delta',b})^2.
\]
The same ergodic argument as in Proposition~\ref{prop:variance}, now applied
to marked point processes indexed by Delaunay triangles rather than by
Delaunay edges, yields
\[
\mathbb{E}
\left[
\left(V_{3,N}^{(W)\prime}\right)^2
\,\middle|\,P_1
\right]
\longrightarrow
\sigma_{V_3}^{2}
\qquad\text{a.s.},
\]
where \(\sigma_{V_3}^{2}\in(0,\infty)\) is deterministic. The finiteness of
\(\sigma_{V_3}^{2}\) follows from the same estimates as for
\(\sigma_{V_2}^{2}\). Indeed, the decay of correlations between distant
increments is still governed by Lemma~\ref{Le:bound:corr}, while the
probability bounds for Delaunay configurations follow from
Lemma~\ref{Le:estimatepN}. The additional factors
\((1-R_\Delta^2)^{-1/2}\) are local functions of the shape of the Delaunay
triangle; their possible singularities near degenerate triangles are
integrable in the Poisson--Delaunay typical-cell distribution for
\(2H\in(0,1)\).

It remains to control the first contraction. As above,
\[
f_{N,3}
=
\frac{1}{\sqrt{|DT_N'|}}
\sum_{\Delta\in DT_N'}
\left(
\eta_{\Delta,1}^{\otimes 2}
+
\eta_{\Delta,2}^{\otimes 2}
\right),
\]
and therefore
\[
\|f_{N,3}\otimes_1 f_{N,3}\|_{\mathfrak H^{\otimes 2}}^{2}
\]
is a normalized quadruple sum of products of four correlations of the form
\[
\mathrm{corr}(Z_{\Delta,a},Z_{\Delta',b}),
\qquad
a,b\in\{1,2\}.
\]
Using the four identities displayed above, each such product is bounded by a
finite sum of products of correlations between normalized edge increments,
multiplied by local geometric factors attached to the involved Delaunay
triangles. The argument of Section~\ref{Sect_Cond_c} applies term by term:
large Delaunay edges are negligible by the exponential tail bound of
Lemma~\ref{Le:estimatepN}, and the contribution of moderate edges is
controlled by the same summability estimates based on the decay
\(\|x-y\|^{2H-2}\). Hence
\[
\|f_{N,3}\otimes_1 f_{N,3}\|_{\mathfrak H^{\otimes 2}}
\longrightarrow 0
\qquad\text{a.s.}
\]

The second-chaos normal approximation criterion therefore gives, conditionally
on \(P_1\),
\[
V_{3,N}^{(W)\prime}
=
I_2(f_{N,3})
\overset{\mathcal D}{\longrightarrow}
\mathcal N(0,\sigma_{V_3}^{2})
\qquad\text{a.s.}
\]
Since the limiting variance is deterministic, the unconditional convergence
follows by dominated convergence of the conditional characteristic
functions. Finally, \(V_{3,N}^{(W)}\) and \(V_{3,N}^{(W)\prime}\) have the
same distribution, and therefore
\[
V_{3,N}^{(W)}
\overset{\mathcal D}{\longrightarrow}
\mathcal N(0,\sigma_{V_3}^{2}).
\]

\section{Technical lemmas}
\label{sec:intermediary_results}

In this section, we establish technical results which are useful to derive Theorem \ref{th:CLTgaussian}.

\subsection{Asymptotic correlations between pairs of normalized increments \label{Sect_Corr_Inc}}

\begin{lemma}
\label{Le:bound:corr}
Let $H\in(0,1)$.
\begin{enumerate}[(i)]
\item Let $d_{1,2}$ and $d_{3,4}$ be fixed and let
$x_{1},x_{2},x_{3},x_{4}$ be such that
\[
\|x_{2}-x_{1}\|=d_{1,2},
\qquad
\|x_{4}-x_{3}\|=d_{3,4}.
\]
Set $d=d_{1,3}:=\|x_{3}-x_{1}\|$. Then, as $d\to\infty$,
\begin{equation*}
\mathrm{corr}
\left(
U_{x_{1},x_{2}}^{(W)},U_{x_{3},x_{4}}^{(W)}
\right)
=
H
(d_{1,2}d_{3,4})^{1-H}
d_{1,3}^{2H-2}
\left(
\cos\beta\cos\theta
-
(1-2H)\sin\beta\sin\theta
\right)
+
o(d_{1,3}^{2H-2}),
\end{equation*}
where
\[
\theta=\mathrm{angle}(\vec u,\overrightarrow{x_{1}x_{2}}),
\qquad
\beta=\mathrm{angle}(\vec u,\overrightarrow{x_{3}x_{4}}),
\]
and where $\vec u$ is a unit vector orthogonal to
$\overrightarrow{x_{3}x_{1}}$ such that
$(\vec u,\overrightarrow{x_{3}x_{1}})$ is positively oriented.

\item Let $\varepsilon\in(0,1/2)$. Then there exist constants
$c>0$ and $d_{0}>0$ such that, for any
$x_{1},x_{2},x_{3},x_{4}\in\mathbf{R}^{2}$ satisfying
\[
0<\|x_{4}-x_{3}\|
\leq
\|x_{2}-x_{1}\|
\leq
\|x_{3}-x_{1}\|^{\varepsilon}
\]
and $\|x_{3}-x_{1}\|\geq d_{0}$, one has
\[
\left|
\mathrm{corr}
\left(
U_{x_{1},x_{2}}^{(W)},U_{x_{3},x_{4}}^{(W)}
\right)
\right|
\leq
c\,
\|x_{2}-x_{1}\|^{2-2H}
\|x_{3}-x_{1}\|^{2H-2}.
\]
\end{enumerate}
\end{lemma}

\begin{prooft}{Lemma \ref{Le:bound:corr}}

Set
\[
a=x_{2}-x_{1},
\qquad
b=x_{4}-x_{3},
\qquad
r=x_{1}-x_{3},
\]
and write
\[
\ell_{1}=\|a\|=d_{1,2},
\qquad
\ell_{2}=\|b\|=d_{3,4},
\qquad
d=\|r\|=d_{1,3}.
\]
From the covariance function \eqref{eq:defcovariance}, we have
\begin{multline*}
\mathrm{cov}
\left(
W(x_{2})-W(x_{1}),
W(x_{4})-W(x_{3})
\right)
\\
=
\frac{\sigma^{2}}{2}
\left(
\|r-b\|^{2H}
-\|r\|^{2H}
-\|r+a-b\|^{2H}
+\|r+a\|^{2H}
\right).
\end{multline*}
Therefore
\begin{equation}
\mathrm{corr}
\left(
U_{x_{1},x_{2}}^{(W)},U_{x_{3},x_{4}}^{(W)}
\right)
=
\frac{\Psi(r,a,b)}
{2(\ell_{1}\ell_{2})^{H}},
\label{eq:corr-Psi}
\end{equation}
where
\[
\Psi(r,a,b)
=
\|r-b\|^{2H}
-\|r\|^{2H}
-\|r+a-b\|^{2H}
+\|r+a\|^{2H}.
\]

Let
\[
F_{r}(z)=\|r+z\|^{2H}.
\]
Then
\[
\Psi(r,a,b)
=
F_{r}(-b)-F_{r}(0)-F_{r}(a-b)+F_{r}(a).
\]
By applying the fundamental theorem of calculus twice, we obtain the exact
identity
\begin{equation}
\Psi(r,a,b)
=
\int_{0}^{1}\int_{0}^{1}
a^{\top}
\nabla^{2}F_{r}(ta-sb)
b
\,\mathrm{d}s\,\mathrm{d}t .
\label{eq:second-order-identity}
\end{equation}
Moreover,
\[
\nabla^{2}F_{r}(z)
=
2H\|r+z\|^{2H-2}I_{2}
+
2H(2H-2)
\|r+z\|^{2H-4}
(r+z)(r+z)^{\top}.
\]

We first prove (i). Since $\ell_{1}$ and $\ell_{2}$ are fixed, uniformly in
$s,t\in[0,1]$,
\[
d^{2-2H}\nabla^{2}F_{r}(ta-sb)
\longrightarrow
2H\left(I_{2}+(2H-2)ee^{\top}\right),
\qquad
e=\frac{r}{\|r\|}.
\]
Choose coordinates such that $e=(0,1)$ and such that the first coordinate
axis is the vector $\vec u$ appearing in the statement. Then
\[
a=\ell_{1}(\cos\theta,\sin\theta),
\qquad
b=\ell_{2}(\cos\beta,\sin\beta).
\]
It follows from \eqref{eq:second-order-identity} that
\begin{align*}
\Psi(r,a,b)
&=
2H d^{2H-2}
\left[
a\cdot b
+
(2H-2)(a\cdot e)(b\cdot e)
\right]
+
o(d^{2H-2})
\\
&=
2H \ell_{1}\ell_{2} d^{2H-2}
\left[
\cos\beta\cos\theta
+
(2H-1)\sin\beta\sin\theta
\right]
+
o(d^{2H-2})
\\
&=
2H \ell_{1}\ell_{2} d^{2H-2}
\left[
\cos\beta\cos\theta
-
(1-2H)\sin\beta\sin\theta
\right]
+
o(d^{2H-2}).
\end{align*}
Combining this expansion with \eqref{eq:corr-Psi} proves (i).

We now prove (ii). Assume that
\[
0<\ell_{2}\leq \ell_{1}\leq d^{\varepsilon},
\qquad
\varepsilon\in(0,1/2).
\]
For $s,t\in[0,1]$, we have
\[
\|ta-sb\|
\leq
\ell_{1}+\ell_{2}
\leq
2d^{\varepsilon}.
\]
Choosing $d_{0}$ large enough, we may ensure that, for all $d\geq d_{0}$,
\[
d-2d^{\varepsilon}\geq \frac{d}{2}.
\]
Hence
\[
\|r+ta-sb\|\geq \frac{d}{2},
\qquad
s,t\in[0,1].
\]
From the explicit expression of the Hessian, there exists a constant
$c>0$, depending only on $2H$, such that
\[
\sup_{s,t\in[0,1]}
\left\|
\nabla^{2}F_{r}(ta-sb)
\right\|
\leq
c d^{2H-2}.
\]
Using \eqref{eq:second-order-identity}, we get
\[
|\Psi(r,a,b)|
\leq
c\,\ell_{1}\ell_{2}\,d^{2H-2}.
\]
Consequently, by \eqref{eq:corr-Psi},
\[
\left|
\mathrm{corr}
\left(
U_{x_{1},x_{2}}^{(W)},U_{x_{3},x_{4}}^{(W)}
\right)
\right|
\leq
c
(\ell_{1}\ell_{2})^{1-H}
d^{2H-2}.
\]
Since $\ell_{2}\leq \ell_{1}$ and $2H<2$,
\[
(\ell_{1}\ell_{2})^{1-H}
\leq
\ell_{1}^{2-2H}.
\]
Therefore
\[
\left|
\mathrm{corr}
\left(
U_{x_{1},x_{2}}^{(W)},U_{x_{3},x_{4}}^{(W)}
\right)
\right|
\leq
c\,\ell_{1}^{2-2H}d^{2H-2},
\]
which is the desired bound.
\end{prooft}

\subsection{Bounds for the density functions of Delaunay neighbors\label%
{Sect_Proof_Le:estmatepN}}

\begin{lemma}
\label{Le:estimatepN}
\begin{enumerate}[(i)]
\item Let \(x_{1},x_{2},x_{3},x_{4}\in \mathbf{R}^{2}\). With the same
notation as in Section \ref{sec:conditional_variance}, for any
\((j,i)\in\mathcal{P}_{2}\),
\[
q_{2,N}^{(j\leftrightarrow i)}
\left(\vec{x}_{\{1:4\}\setminus\{j\}}\right)
\leq
\left(\pi N R^{2}+4\right)
\exp\left\{
-\frac{\pi}{4}NR^{2}
\right\},
\]
where
\[
R=
\max\left\{
\|x_{2}-x_{1}\|,
\|x_{\{3,4\}\setminus\{j\}}-x_{i}\|
\right\}.
\]

\item Let \(x_{1},x_{2},x_{3},x_{4}\in \mathbf{R}^{2}\), and let
\(p_{2,N}(x_{1},x_{2},x_{3},x_{4})\) be as in Eq.~\eqref{eq:defp2N}.
Assume that
\[
\|x_{4}-x_{3}\|\leq \|x_{2}-x_{1}\|.
\]
Then
\[
p_{2,N}(x_{1},x_{2},x_{3},x_{4})
\leq
\left(
\pi N\|x_{2}-x_{1}\|^{2}+4
\right)
\exp\left\{
-\frac{\pi}{4}N\|x_{2}-x_{1}\|^{2}
\right\}.
\]
\end{enumerate}
\end{lemma}

\begin{prooft}{Lemma \ref{Le:estimatepN}}

We first prove \textit{(ii)}. Set
\[
L=\|x_{2}-x_{1}\|.
\]
Since \(\|x_{4}-x_{3}\|\leq L\), on the event defining
\(p_{2,N}(x_{1},x_{2},x_{3},x_{4})\), the edge \([x_{1},x_{2}]\) is a
Delaunay edge in the triangulation generated by
\(P_N\cup\{x_{1},x_{2},x_{3},x_{4}\}\). Removing the additional fixed points
can only make the empty-circle condition easier to satisfy. Hence
\[
p_{2,N}(x_{1},x_{2},x_{3},x_{4})
\leq
\mathbb{P}
\left[
x_{1}\sim x_{2}
\text{ in }
\operatorname{Del}(P_N\cup\{x_{1},x_{2}\})
\right].
\]
If \(x_{1}\) and \(x_{2}\) are Delaunay neighbors in
\(\operatorname{Del}(P_N\cup\{x_{1},x_{2}\})\), then, almost surely, there
exists \(y\in P_N\) such that
\[
\Delta(x_{1},x_{2},y)\in
\operatorname{Del}(P_N\cup\{x_{1},x_{2}\}).
\]
Consequently, by a union bound and the Slivnyak--Mecke formula,
\begin{align*}
p_{2,N}(x_{1},x_{2},x_{3},x_{4})
&\leq
\mathbb{E}
\left[
\sum_{y\in P_N}
\mathbb{I}
\left[
P_N\cap B(x_{1},x_{2},y)=\emptyset
\right]
\right]
\\
&=
N\int_{\mathbf{R}^{2}}
\exp\left\{
-N a(B(x_{1},x_{2},y))
\right\}
\,\mathrm{d}y,
\end{align*}
where \(B(x_{1},x_{2},y)\) denotes the circumdisk passing through
\(x_{1},x_{2},y\). The collinear case is irrelevant since it has Lebesgue
measure zero in the integral.

The radius of the circumdisk \(B(x_{1},x_{2},y)\) is at least
\[
\frac12
\max\left\{
\|x_{2}-x_{1}\|,
\|y-x_{1}\|
\right\}.
\]
Therefore
\[
a(B(x_{1},x_{2},y))
\geq
\frac{\pi}{4}
\max\left\{
L,\|y-x_{1}\|
\right\}^{2},
\]
and thus
\[
p_{2,N}(x_{1},x_{2},x_{3},x_{4})
\leq
N
\int_{\mathbf{R}^{2}}
\exp\left\{
-\frac{\pi}{4}N
\max\left\{
L,\|y-x_{1}\|
\right\}^{2}
\right\}
\,\mathrm{d}y.
\]
Splitting the integral according to whether \(\|y-x_{1}\|\leq L\) or
\(\|y-x_{1}\|>L\), we obtain
\begin{align*}
p_{2,N}(x_{1},x_{2},x_{3},x_{4})
&\leq
N e^{-\frac{\pi}{4}NL^{2}}
\int_{\mathbf{R}^{2}}
\mathbb{I}\left[\|y-x_{1}\|\leq L\right]
\,\mathrm{d}y
\\
&\quad
+
N
\int_{\mathbf{R}^{2}}
e^{-\frac{\pi}{4}N\|y-x_{1}\|^{2}}
\mathbb{I}\left[\|y-x_{1}\|>L\right]
\,\mathrm{d}y
\\
&=
\pi N L^{2}e^{-\frac{\pi}{4}NL^{2}}
+
4e^{-\frac{\pi}{4}NL^{2}}.
\end{align*}
This proves
\[
p_{2,N}(x_{1},x_{2},x_{3},x_{4})
\leq
\left(
\pi N L^{2}+4
\right)
e^{-\frac{\pi}{4}NL^{2}},
\]
which is the desired bound.

We now prove \textit{(i)}. On the event defining
\(q_{2,N}^{(j\leftrightarrow i)}\), two Delaunay edges share one endpoint.
Let \(R\) be the larger of their two lengths. The event implies, in
particular, that the edge of length \(R\) is a Delaunay edge in the
triangulation generated by the Poisson point process together with the fixed
points under consideration. Removing the remaining fixed point can only
preserve the existence of an empty circumdisk for this edge. Applying the
single-edge bound just proved, with \(L=R\), gives
\[
q_{2,N}^{(j\leftrightarrow i)}
\left(\vec{x}_{\{1:4\}\setminus\{j\}}\right)
\leq
\left(\pi N R^{2}+4\right)
\exp\left\{
-\frac{\pi}{4}NR^{2}
\right\}.
\]
This proves \textit{(i)}.
\end{prooft}

\subsection{Bounds for some exponential moments of a uniform distribution
over a disc}

Let \(N>0\) and let \(R\) be a positive random variable with distribution
function
\begin{equation}
\mathbb{P}[R\leq r]
=
\begin{cases}
0, & r<0,\\[2mm]
\dfrac{r^{2}}{N}, & 0\leq r\leq \sqrt{N},\\[2mm]
1, & r>\sqrt{N}.
\end{cases}
\label{Eq_Dist_R}
\end{equation}

\begin{lemma}
\label{Le_Bound_Exp_R}
Let \(0<H<1/2\) and \(d_{0}>0\). There exist two constants
\(c_{1}\) and \(c_{2}\), depending only on \(2H\) and \(d_{0}\), such
that, for large \(N\),
\[
\mathbb{E}
\left[
\exp\left(R^{2H-2}\mathbb{I}[R\geq d_{0}]\right)
\right]
\leq
1+c_{1}(\sqrt{N})^{2H-2}+c_{2}N^{-1}.
\]
\end{lemma}

\begin{prooft}{Lemma \ref{Le_Bound_Exp_R}}
For \(N>d_{0}^{2}\), set \(a_{N}=d_{0}/\sqrt{N}\). Since
\(R/\sqrt{N}\) has density \(2r\mathbb{I}_{[0,1]}(r)\), we have
\begin{align*}
\mathbb{E}
\left[
\exp\left(R^{2H-2}\mathbb{I}[R\geq d_{0}]\right)
\right]
&=
\int_{0}^{a_N}2r\,\mathrm{d}r
+
\int_{a_N}^{1}
\exp\left((r\sqrt{N})^{2H-2}\right)
2r\,\mathrm{d}r
\\
&=
1+
2\int_{a_N}^{1}
\left[
\exp\left((r\sqrt{N})^{2H-2}\right)-1
\right]r\,\mathrm{d}r .
\end{align*}
Using the expansion \(e^x-1\leq x+\sum_{k=2}^{\infty}x^k/k!\) for
\(x\geq0\), we get
\begin{align*}
\mathbb{E}
\left[
\exp\left(R^{2H-2}\mathbb{I}[R\geq d_{0}]\right)
\right]
&\leq
1+
2\int_{a_N}^{1}
(r\sqrt{N})^{2H-2}r\,\mathrm{d}r
\\
&\quad
+
2\int_{a_N}^{1}
\sum_{k=2}^{\infty}
\frac{(r\sqrt{N})^{k(2H-2)}}{k!}
r\,\mathrm{d}r .
\end{align*}
For the first integral,
\[
2\int_{a_N}^{1}
(r\sqrt{N})^{2H-2}r\,\mathrm{d}r
=
2(\sqrt{N})^{2H-2}
\int_{a_N}^{1}r^{2H-1}\,\mathrm{d}r
\leq
\frac{2}{2H}(\sqrt{N})^{2H-2}.
\]

For the second integral, since \(r\sqrt{N}\geq d_{0}\) on
\([a_N,1]\) and \(2H-2<0\), we have
\[
\sum_{k=2}^{\infty}
\frac{(r\sqrt{N})^{k(2H-2)}}{k!}
\leq
(r\sqrt{N})^{2(2H-2)}
\exp\left((r\sqrt{N})^{2H-2}\right)
\leq
(r\sqrt{N})^{2(2H-2)}
\exp\left(d_{0}^{2H-2}\right).
\]
Therefore
\begin{align*}
2\int_{a_N}^{1}
\sum_{k=2}^{\infty}
\frac{(r\sqrt{N})^{k(2H-2)}}{k!}
r\,\mathrm{d}r
&\leq
2e^{d_{0}^{2H-2}}
N^{2H-2}
\int_{a_N}^{1}
r^{4H-3}\,\mathrm{d}r
\\
&=
\frac{e^{d_{0}^{2H-2}}}{1-2H}
N^{2H-2}
\left(a_N^{4H-2}-1\right)
\\
&\leq
\frac{e^{d_{0}^{2H-2}}}{1-2H}
d_{0}^{2(2H-1)}
N^{-1}.
\end{align*}
Combining the two estimates yields
\[
\mathbb{E}
\left[
\exp\left(R^{2H-2}\mathbb{I}[R\geq d_{0}]\right)
\right]
\leq
1+
\frac{1}{H}(\sqrt{N})^{2H-2}
+
\frac{e^{d_{0}^{2H-2}}}{1-2H}
d_{0}^{2(2H-1)}
N^{-1}.
\]
Thus the result holds with
\[
c_{1}=\frac{1}{H},
\qquad
c_{2}=
\frac{e^{d_{0}^{2H-2}}}{1-2H}
d_{0}^{2(2H-1)}.
\]
\end{prooft}

 \end{appendices}

%

\end{document}